\documentclass[article]{IEEEtran}
\IEEEoverridecommandlockouts  

\usepackage{graphics} 
\usepackage{epsfig} 
\usepackage{times} 
\usepackage{amsmath} 
\usepackage{amssymb}  
\usepackage{dsfont}
\usepackage{mathtools}
\usepackage[utf8]{inputenc}
\usepackage{graphicx}
\usepackage{newtxtext, newtxmath}
\usepackage{color,soul}
\usepackage{comment}
\allowdisplaybreaks

\usepackage{enumerate}
\newcommand{\R}{\mathbb{R}}

\newcommand{\M}{\mathbb{M}}
\newcommand{\W}{\mathbb{W}}

\newcommand{\norm}[1]{\|#1\|}
\let\phi\varphi
\newcommand{\abs}[1]{|#1|}
\newcommand{\X}{\mathbb{X}}
\newcommand{\U}{\mathbb{U}}
\newcommand{\Z}{\mathbb{Z}}

\DeclareMathOperator*{\argmin}{arg\,min}

\newtheorem{definition}{Definition}[section]
\newtheorem{theorem}{Theorem}[section]
\newtheorem{remark}{Remark}[section]

\newtheorem{lemma}{Lemma}[section]
\newtheorem{assumption}{Assumption}[section]

\allowdisplaybreaks

\def\qed{\hfill $\diamond$}

\usepackage[colorlinks=true,breaklinks=true,bookmarks=true,urlcolor=blue,citecolor=blue,linkcolor=blue,bookmarksopen=false,draft=false]{hyperref}

\begin{document}
\title{
Reinforcement Learning for Jointly Optimal Coding and Control Policies for a Controlled Markovian System over a Communication Channel
  \thanks{This work was supported by the Natural Sciences and Engineering Research Council of Canada. The first and third authors are with the Department of Mathematics and Statistics, Queen's University, Kingston ON, Canada (email: \{19eah7,yuksel\}@queensu.ca). The second author is with the Department of Information Technology and Electrical Engineering, ETH Z{\"u}rich, Switzerland (email: lcregg@ethz.ch).}}

\author{Evelyn Hubbard, Liam Cregg, Serdar Y\"uksel}
\maketitle

\begin{abstract}
We study the problem of joint optimization involving coding and control policies for a controlled Markovian system over a finite-rate noiseless communication channel. While structural results on the optimal encoding and control have been obtained in the literature, their implementation has been prohibitive in general, except for linear models. To this end, we first develop existence, regularity, and structural properties on optimal policies, followed by rigorous approximations and reinforcement learning results. Notably, we establish near optimality of finite model approximations obtained via predictor (conditional probability on the next state realization given the information at the controller) quantization as well as sliding finite window approximations, and the convergence of a reinforcement learning algorithm to near optimality. A detailed comparison of the approximation schemes and their performance is presented. To our knowledge, this is the first rigorous reinforcement learning study on jointly optimal coding and networked control with finite rate channels. 
\end{abstract}

\section{Introduction}
Networked control systems involve stochastic control systems where a communication channel is present between different stations, such as sensors, actuators and controllers. In this context there exists a comprehensive literature on stabilization and optimization of such systems under various information constraints, see e.g. \cite{franceschetti2014elements,Hespanha,MatveevSavkin,YukselBasarBook24,kawan2013invariance} for extensive reviews.

Networked control theory requires interdisciplinary methods to arrive at optimality and structural results on optimal coders, decoders and controllers. While analytically and architecturally very useful, the structural results on optimal policies often lead to optimization problems in uncountable spaces and hence lead to computationally challenging problems. Furthermore, despite the presence of such structural results in the literature as we will review, learning theoretic and rigorously justified near optimality results have been limited. The goal of this paper is to present reinforcement learning algorithms which are provably near optimal for the jointly optimal coding and control policies in a system which is connected to a controller over a finite rate noiseless channel.

\subsection{Problem Setup} \label{problemsetup}
We consider a networked control problem where a controlled Markov source observed over a noiseless communication channel is controlled using data obtained from this channel. The controlled Markov source is updated at each time step with the following dynamics. For $t\geq0,$
\begin{align} \label{system}
    x_{t+1} = f(x_t,u_t,w_t),
\end{align}
where $x_t$ takes values from a finite\footnote{The case with continuous state spaces is discussed in Section \ref{contSpaceSec}.} state space $\X$, $u_t$ takes values from a finite action space $\mathbb{U}$, and $w_t$ is some independent and identically distributed (i.i.d.) noise process. The initial state $x_0$ is a random variable with distribution $\pi_0$. 

\begin{figure}[h] \label{sys}
\begin{center}
\includegraphics[width = 9cm]{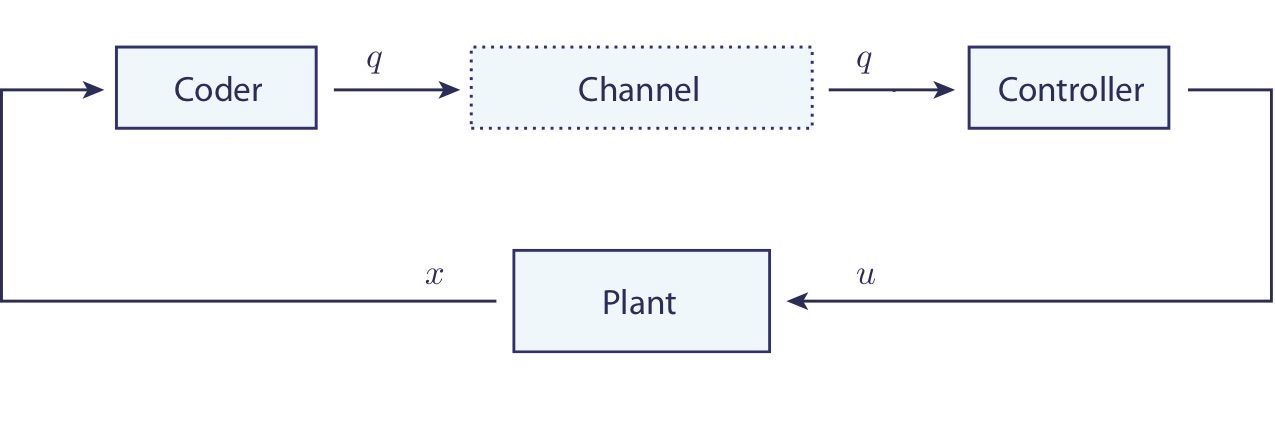}
\caption{Control driven Markov process over noiseless channel 
\label{loop}}
\end{center}
\end{figure}

At each time stage $t$, $x_t$ is encoded causally over a noiseless channel. The controller receives information from the communication channel and chooses a control action $u_t$ which is then transmitted to the plant. The encoder and controller use coding and control policies, respectively, to compute their outputs. This system is depicted in Figure \ref{loop}.
\\\indent A \textit{coding policy} $\gamma^{e} = \{\gamma^e_t, t \ge 0\}$ is a sequence of functions which generate quantization outputs, $q_t$, as a measurable function of the encoder's information at time $t$:
\begin{equation}
    I_t^e = \{{x_{[0,t]},q_{[0,t-1]}}\}
\end{equation}
\begin{equation}
    \gamma^e_t:I_t^e\mapsto q_t \in \M
\end{equation}
where $\M$ is a finite quantization output alphabet, so that $q_t = \gamma^e_t(I_t^e)$ for $t \in \mathbb{Z}_+$.


A \textit{control policy} $\gamma^c = \{\gamma^c_t,t \ge 0\}$ is a sequence of functions which generate control actions, $u_t$, as a measurable function of the information available at the controller at time $t$:
\begin{equation}
    I_t^c = \{{q_{[0,t]}, u_{[0,t-1]}}\}
\end{equation}
\begin{equation}
    \gamma^c_t:I_t^c \mapsto u_t \in \mathbb{U} 
\end{equation}
where $\mathbb{U}$ is the finite action space, so that $u_t = \gamma^c_t(I_t^c)$ for all $t \in \mathbb{Z}_+$.
\begin{definition}
With the coding and control policies as given above, we define the class of admissible \textit{joint coding-control policies} as
    \begin{equation}
        \Gamma_A \coloneqq \{\overline{\gamma} = \{(\gamma^e_t,\gamma^c_t), t\geq0\} \}.
    \end{equation}
\end{definition}
The networked control problem involves a joint optimization of coding and control policies\footnote{A noisy channel variation of this problem with feedback as depicted in Figure \ref{LLL0001} is studied in Section \ref{noisyChannelExt}.}. The optimization objective for a finite time horizon is defined as follows for some $N \ge 0$:
\begin{equation} \label{objective}
    J^N(\pi_0) = \inf_{\overline\gamma\in\Gamma_A} J^N(\pi_0,\overline{\gamma}) \coloneqq  \inf_{\overline{\gamma}\in\Gamma_A} E_{\pi_0}^{\overline{\gamma}}\left[\sum^{N-1}_{k=0} c(x_k,u_k)\right],
\end{equation} 
where $c(x_t,u_t)$ is a predefined cost function at each state-action pair. Here we use $E^{\overline{\gamma}}_{\pi_0}$ to denote the expectation on $(x_t, u_t)_{t \ge 0}$ given initial distribution $\pi_0$ and joint coding-control policy $\overline{\gamma}$. Our focus will be on an infinite horizon discounted cost criterion to be presented further below (\ref{infinitediscountcost}). 

    \subsection{Literature review and contributions}
Before proceeding further, we acknowledge and review several related studies, primarily in the control-free domain. In the control-free setup, related papers on real-time coding include the following: \cite{NeuhoffGilbert} established that the optimal causal encoder minimizing the data rate subject to a distortion constraint for an  i.i.d. sequence is memoryless. If the source is $k$th-order Markov, then the optimal causal fixed-rate coder minimizing any measurable distortion uses only the last $k$ source symbols, together with the current state at the receiver's memory \cite{Witsenhausen}. Reference \cite{WalrandVaraiya} considered the optimal causal coding problem of finite-state Markov sources over noisy channels with feedback. \cite{Teneketzis} and \cite{MahTen09} considered optimal causal coding of Markov sources over noisy channels without feedback. \cite{MahajanTeneketzisJSAC} considered the optimal causal coding over a noisy channel with noisy feedback. Reference \cite{LinderZamir} considered the causal coding of more general sources, stationary sources, under a high-rate assumption. An earlier reference on quantizer design is \cite{Curry}. Relevant discussions on optimal quantization, randomized decisions, and optimal quantizer design can be found in \cite{GaborGyorfi} and \cite{YukselOptimizationofChannels}. \cite{BorkarMitterTatikonda} has studied a related problem of coding of a partially observed Markov source.  \cite{NayyarTeneketzis} considered within a multi-terminal setup decentralized coding of correlated sources when the encoders observe conditionally independent messages given a finitely valued random variable, and obtained separation results for the optimal encoders. 

Further related studies include sequential decentralized hypothesis testing problems \cite{VeeravalliBasarPoor} and multi-access communications with feedback \cite{AchilleasAllerton09}. A detailed review is available in \cite[Chapter 15]{YukselBasarBook24}. 

Quantizer design for the linear quadratic Gaussian setup has been studied by many \cite{LewisTou,Curry,Fischer82,BaoMikaelKalle,MatveevSavkin04,TatikondaSahaiMitter,NairFagnaniSurvey,BorkarMitterLQG,FuTAC2012,rabi2016separated,YukselLQGQuantizationTAC,yuksel2019note}. Information theoretic methods lead to further methods, though with operational restrictions compared with zero-delay schemes: \cite{banbas89,TatikondaSahaiMitter,charalambous2014nonanticipative,tanaka2015semidefinite,tanaka2016semidefinite,tanaka2018lqg,stavrou2018zero,derpich2012improved,stavrou2018optimal,stavrou2021asymptotic,kostina2019rate,cuvelier2022time}. 

Finally, \cite{creggZeroDelayNoiseless} and \cite{cregg2024reinforcementlearningoptimaltransmission} studied approximation and reinforcement learning techniques which are closely relevant to those studied in this paper, in the control-free case. However, the controlled setup requires a different analysis since the decoder/controller structure can be taken a priori to be a Bayesian minimizer in the control-free setup (and which can be computed off-line), which is not true in the controlled case. In the controlled case, a more delicate and significantly more general MDP formulation is needed with associated continuity and regularity analysis, and, in particular, the learning analysis requires that the controller also be optimized unlike the control-free case. 

{\bf Contributions.} 
\begin{itemize}
\item[(i)] While structural results on the optimal encoding and control have been obtained in the literature especially for finite horizon problems, their implementation has been prohibitive in general, except for linear models. In this paper, we develop several regularity and structural properties leading to optimality results for jointly optimal coding and control.
\item[(ii)] We establish regularity properties such as weak Feller continuity (Theorem \ref{wf}), and filter/predictor stability (Theorem \ref{loss bound}). Building on these, we show that the dynamic programming recursions are well defined and an optimal solution exists under both finite and infinite horizon discounted optimization criteria (Theorem \ref{best performance class}). These regularity results lead to rigorous near optimal approximation results where approximations are obtained either with quantized approximations or are obtained by finite window control policies, leading to near optimal solutions under complementary conditions (Theorem \ref{special case quantization performance} and Theorem \ref{T14}, respectively). Accordingly, both existence and near optimality results are new contributions to the literature, to our knowledge.
\item[(iii)] Finally, corresponding to both of the approximation methods, we obtain rigorous reinforcement learning results obtained via associated quantized Q-learning algorithms and its convergence to near optimality: these involve reinforcement learning where the variables used in the algorithm are either the quantization of the predictor variables or are sliding finite window of the most recent measurement and action variables. Both of the algorithms are shown to converge to near optimality under complementary conditions (Theorem \ref{theorem:quantizedQlearning_nearoptimal} and Theorem \ref{theorem:finite_memory_Qlearning}). We present the comparative benefits of either of the methods in Section \ref{comparisonModelsSec}.  
\end{itemize}

\subsection{Notation} We use both uppercase and lowercase letters to denote random variables; the distinction between a random variable and its realization will either be clear from context or, when important, explicitly identified. For a given sequence $(x_t)_{t \ge 0}$, we denote a contiguous subset $(x_n, x_{n+1}, \ldots, x_m)$ by $x_{[n,m]}$. We use $P(\cdot)$ and $E[\cdot]$ to denote probability measures and expectations of a given random variable, respectively. When these depend on certain parameters, we include these parameters in the superscript and/or subscript. When integrating against a given probability measure, say $P$, we will make explicit the random variable that the measure corresponds to; that is, we write $\int g(x_t) P(dx_t) \coloneqq \int g dP$, when $P$ is the probability measure corresponding to the random variable $x_t$ (thus the $x_t$ in the integral should be interpreted as a realization; as mentioned earlier, $x_t$ may be either a random variable or its realization depending on the context). We may use this integral notation even when the relevant random variable is finite (and thus the expectation would be a sum); in this case $P(dx_t)$ is the appropriate counting measure. For a given Polish (that is, complete, separable, and metric) space $\mathbb{X}$, we denote the set of probability measures over $\mathbb{X}$ by $\mathcal{P}(\mathbb{X})$, and we denote the Borel $\sigma$-algebra over $\mathbb{X}$ by $\mathcal{B}(\mathbb{X})$. When the time index is not important, we write Markov transition probabilities $P(dx_{t+1} | x_t, u_t)$ as $P(dx' | x,u)$. Finally, for a metric space $\mathbb{X}$, $C_b(\mathbb{X})$ denotes the set of all continuous and bounded functions from $\mathbb{X}$ to $\mathbb{R}$.
\begin{table}[htbp]
\centering
\caption{Summary of Notation}
\renewcommand{\arraystretch}{1.2}
\begin{tabular}{ll}
\textbf{Symbol} & \textbf{Description} \\ 
$\X$, $\mathbb{U}$ & Finite state and action spaces \\ 
$x_t \in \X$ & System (plant) state at time $t$ \\ 
$u_t \in \mathbb{U}$ & Control action at time $t$ \\ 
$w_t$ & i.i.d. noise process. \\ 
$\pi_t \in \mathcal{P}(\X)$ & Predictor at time $t$ \\ 
$\bar{\pi}_t \in \mathcal{P}(\X)$ & Filter at time $t$\\
$Q_t : \X \to \M$ & Quantizer chosen at time $t$ \\ 
$q_t \in \M$ & Quantization output at time $t$ \\ 
$\eta_t : \mathcal{Q} \times \M \to \mathbb{U}$ & Control mapping at time $t$ \\ 
$\gamma^e_t$, $\gamma^c_t$ & Encoding and control policies at time $t$ \\ 
$\Gamma_A$ & Set of admissible coding-control policies \\ 
$\Gamma_{C\!-\!P}$ & Class of predictor-structured policies \\ 
$P(dx'|x,u)$ & Transition kernel of the Markov source \\ 
$c(x,u)$ & Stage cost for state-action pair $(x,u)$ \\ 
$P(d\pi'|\pi,Q,\eta)$ & Transition kernel of the predictor process \\ 
$\tilde{c}(\pi, Q, \eta)$ & Equivalent cost in the predictor-structured MDP \\ 
$J^N(\pi_0,\gamma)$ & Finite-horizon cost under policy $\gamma$ and initial $\pi_0$ \\ 
$J_\beta(\pi_0,\gamma)$ & Infinite-horizon discounted cost under policy $\gamma$ \\ 
$J_\beta(\pi_0)$ & Infinite-horizon discounted cost under optimal policy\\
$\beta \in (0,1)$ & Discount factor \\ 
$F(\pi_t,Q_t,\eta_t,q_t)$ & Predictor update\\
$\mathcal{Q}$, $\mathcal{H}$ & Spaces of quantizers and control mappings, respectively \\ 
$\hat{\pi}_t$ & Approximated predictor at time $t$\\

\end{tabular}
\end{table}
\section{Structural Results and Equivalent MDP Formulation}\label{section:MDP_Formulation}
\subsection{Optimal Structure of Encoding and Control Policies}\label{optimalstructure}
Define the predictor and filter sequences, respectively as
\begin{equation} \label{predictor}
    \pi_t (\cdot) := P^{\overline{\gamma}}(x_t=\cdot\,|q_{[0,t-1]})
\end{equation}
\begin{equation}\label{filter_2}
    \overline{\pi}_t(\cdot) = P^{\overline{\gamma}}(x_t=\cdot\,|q_{[0,t]}).
\end{equation}

 A celebrated structural result of Witsenhausen on the structure of optimal encoders \cite{Witsenhausen} is extended in \cite{YukselLQGQuantizationTAC} and \cite{YukselBasarBook24} to models driven by control:
\begin{theorem} \cite{Witsenhausen}, \cite{YukselLQGQuantizationTAC}, \cite[Theorem 15.3.6]{YukselBasarBook24}\\ 
    For a system with dynamics \eqref{system}, and optimization problem \eqref{objective}, any joint coding-control policy (with a given control policy) can be replaced, without loss in performance, by one which uses $x_t$ and $q_{[0,t-1]}$ (while not altering the control policy), at each $t\geq0$. That is,
    \begin{align*}
        q_t &= \gamma^e_t(q_{[0,t-1]}, x_t),
    \end{align*} 
    \end{theorem}

Walrand and Varaiya further refined Witsenhausen's structural results in \cite{WalrandVaraiya,WalrandVaraiyaControl}, which include extensions to controlled setups (see also \cite{YukselLQGQuantizationTAC} and \cite[Theorem 15.3.6]{YukselBasarBook24}). The version below follows from  \cite[Theorem 15.3.6]{YukselBasarBook24} and its proof.

\begin{theorem}\cite[Theorem 3.2]{WalrandVaraiyaControl}\cite[Theorem 15.3.6]{YukselBasarBook24}\label{Walrand and Varaiya}\\ 
     For a system with dynamics \eqref{system} and optimization problem \eqref{objective}, for any joint coding-control policy, the coding policy can be replaced, without loss in performance, by one which uses only $\pi_t$ and $x_t$, at each $t\geq0$. Furthermore, the control policy can be replaced by one which uses only $\overline{\pi}_t$. That is, without any loss in performance, we have:
     \begin{align}
         q_t & = \gamma^e_t(\pi_t, x_t).\label{eq:WV_coding}
     \end{align}
     
     \begin{equation}
         u_t = \gamma^c_t(\overline{\pi}_t).\label{eq:WV_control}
     \end{equation}
\end{theorem}
It will be convenient for our setup to view \eqref{eq:WV_coding} and \eqref{eq:WV_control} not as functions mapping $(\pi_t,x_t) \mapsto q_t$ and $\overline{\pi}_t \mapsto u_t$, respectively, but instead as functions from $\pi_t$ to an appropriate space of \emph{functions} which are then used to select $q_t$ and $u_t$. First, define the space of \emph{quantizers} as:
$$ \mathcal{Q} \coloneqq \{ Q : \X \to \M \} .$$
Then we can view a coding policy $\gamma^e_t$ in \eqref{eq:WV_coding} equivalently as first choosing a quantizer $Q_t$ using $\pi_t$, then applying this quantizer to $x_t$ to obtain $q_t$. That is,
\begin{align}
    Q_t & = \gamma^e_t(\pi_t) & q_t & = Q_t(x_t).\label{eq:gamma^e}
\end{align}
Note that we can compute $\overline{\pi}_t$ given $(\pi_t, Q_t,q_t)$ using Bayes' formula. So we can similarly view a control policy $\gamma^c_t$ in \eqref{eq:WV_control} as using $\pi_t$ to choose a function $\eta_t : \mathcal{Q} \times \M \to \mathbb{U}$, then applying this function to $(Q_t, q_t)$ to obtain $u_t$. That is,
\begin{align}
    \eta_t & = \gamma^c_t(\pi_t) & u_t & = \eta_t(Q_t, q_t).\label{eq:gamma^c}
\end{align}
We note that no measurability issues arise in this finite setup in the construction of the equivalent policies.

In view of the above, we will define the following class of joint coder-control policies. Let $\mathcal{P}(\X)$ be 
endowed with the weak convergence topology.
\begin{definition}
    A joint coder-controller policy $\overline{\gamma} = \{\gamma^e,\gamma^c\}$, has a \textit{Controlled-Predictor-Structure} if for each $t\geq 0$, $\gamma_t^e$ and $\gamma_t^c$ have the structural forms in \eqref{eq:gamma^e} and \eqref{eq:gamma^c}. That is:
\begin{equation}
    \gamma^{e}_t:\mathcal{P}(\X) \ni \pi_t\mapsto Q_t \in {\cal Q}:=\{Q : \mathbb{X} \to \M\},\label{pi_to_Q}
\end{equation}
and
\begin{equation} 
    \gamma^{c}_t: {\cal P}(\X) \ni \pi_t\mapsto \eta_t \in \mathcal{H} \coloneqq \left\{ \eta : \mathcal{Q}\times\mathcal{M}\to\mathbb{U} \right\}.\label{pi_to_eta}
\end{equation}
We will denote this set of policies as $\Gamma_{C-P}$.
\end{definition}

\begin{figure}[h] \label{reformulated_loop}
\begin{center}
\includegraphics[width=9cm]{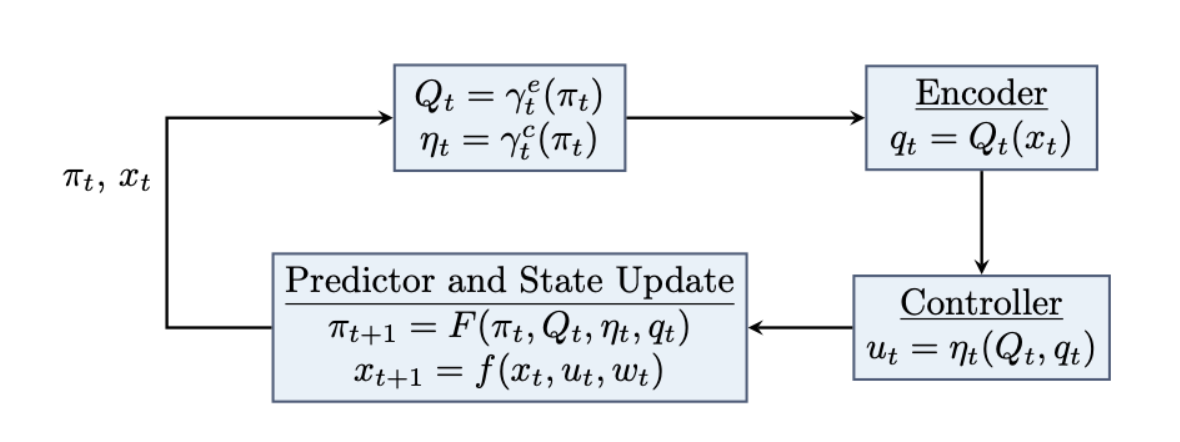}
\caption{Block diagram for the reformulated system: $x_t$ is the
  source sample, $\pi_t$ the predictor on the state, $Q_t$ and $\eta_t$ are quantizer and control maps, respectively chosen by encoder and controller policies, $q_t$ is the encoded symbol transmitted over the
  noiseless channel, and $u_t$ is the control action applied to the system. State and predictor updates are computed and fed back into the loop.}
\end{center}
\end{figure}

Accordingly, a policy in $\Gamma_{C-P}$ will map $\pi_t$ to $(Q_t, \eta_t)$ at time $t$. These are used to generate channel inputs and actions via $q_t = Q_t(x_t)$ and $u_t=\eta_t(Q_t,q_t)$. Thus, the effective {\it state}, as will be justified later in Theorem \ref{costeqiv and blackwell}, is $\pi_t$ and the {\it action} is $(Q_t, \eta_t)$. Under a policy which satisfies the structure given in Theorem \ref{Walrand and Varaiya}, via a Bayesian update, we have the following predictor update equation:
    \begin{align} \label{predictor update}
    &\pi_{t+1}(x_{t+1}) \\
    &= \frac{\sum_{x_t}\sum_{u_t} P(x_{t+1}|x_t,u_t)P(q_t|\pi_t,x_t)P(u_t|q_t,\pi_t)\pi_t(x_t)}{\sum_{x_{t+1}}\sum_{x_t}\sum_{u_t} P(x_{t+1}|x_t,u_t)P(q_t|\pi_t,x_t)P(u_t|q_t,\pi_t)\pi_t(x_t)}\nonumber
    \end{align}
Consider
\begin{equation*}
Q_t^{-1}(q_t) := \{x_t\in \mathbb{X} : Q_t(x_t)=q_t\}
\end{equation*}
Then, we can write the update equation as
\begin{align*}
\pi_{t+1}(x_{t+1})&= \frac{1}{\pi_t(Q_t^{-1}(q_t))} \sum_{Q_t^{-1}(q_t)} P(x_{t+1}|x_t,\eta_t(Q_t,q_t))\pi_t(x_t) \nonumber\\
&:= F(\pi_t,Q_t,\eta_t,q_t)
\end{align*}
We note that this update equation implies that the controller receives partial information (an encoded message) from the communication channel, and uses it to update the predictor of the state accordingly. The update step is effectively the decoding of the message and the controller inherently functions as both a decoder and a decision-maker.\\
The reformulated system is depicted in Figure \ref{reformulated_loop}.
\subsection{Predictor Structured Controlled Markov Problem} \label{Predictor Structured Markov Problem}
The structures for policies in $\Gamma_{C-P}$ motivate the following theorem, extended from \cite[Theorem 15.4.1]{YukselBasarBook24}. Throughout, we assume that we are using such a policy (which is without loss of optimality by Theorem \ref{Walrand and Varaiya}). 
\begin{theorem}\label{Markov property}
    $(\pi_t, (Q_t, \eta_t))$ is a controlled Markov chain with $\pi_t$ (defined on $\mathcal{P}(\X)$) as the state and $(Q_t,\eta_t)$ (defined on $\mathcal{Q}\times \mathcal{H}$) as the control action.
\end{theorem}
\textbf{Proof.}
Let $D \in \cal B({\cal P}(\mathbb{X}))$. Then, under a policy in $\Gamma_{C-P}$,
\begin{align*}
    P&(\pi_{t+1} \in D|\pi_s,Q_s,\eta_s,s\leq t)\\
    &= \sum_{q_t\in \M}P(\pi_{t+1}\in D,q_t|\pi_s,Q_s,\eta_s,s\leq t)\\ 
    &= \sum_{q_t\in \M}P(\pi_{t+1}\in D|q_t,\pi_s,Q_s,\eta_s,s\leq t)P(q_t|\pi_s,Q_s,\eta_s,s\leq t)\\ 
    &=\sum_{q_t\in \M} P(F(\pi_t,Q_t,\eta_t,q_t)\in D|\pi_t,Q_t,\eta_t,q_t)P(q_t|\pi_t,Q_t,\eta_t) \\ 
    &= P(F(\pi_t,Q_t,\eta_t,q_t)\in D|\pi_t,Q_t,\eta_t)= P(\pi_{t+1}\in D|\pi_t,Q_t,\eta_t)
\end{align*}
    The third equality is due to the fact that:
\begin{align*}
    P(&q_t|\pi_s,Q_s,\eta_s,s\leq t)\\
    &=\sum_{x_t\in\X} P(q_t,x_t|\pi_s,Q_s,\eta_s,s\leq t)\\ 
    &= \sum_{x_t\in\X}P(q_t|x_t,\pi_s,Q_s,\eta_s,s\leq t)P(x_t|\pi_s,Q_s,\eta_s,s\leq t)\\
    &= \sum_{x_t\in\X} P(q_t|x_t,Q_t)\tilde{P}(x_t|\pi_t,Q_t,\eta_t)
    = P(q_t|\pi_t,Q_t,\eta_t),
\end{align*}
where we note that $q_t$ is determined by $x_t$ and $Q_t$ and since under $\overline{\gamma} \in \Gamma_{C-P}$, conditioning on $\pi_t$ implicitly conditions on $Q_{[0,t-1]}$ and $\eta_{[0,t-1]}$. \qed

Thus we have an MDP with transition probability $P(d\pi_{t+1} | \pi_t, Q_t, \eta_t)$. In order to apply approximation and learning algorithms to this redefined MDP, we build on \cite[Chapter 4]{SaLiYuSpringer} and show that that the transition kernel is weakly continuous \cite{YukLinZeroDelay} (that is, the kernel is weak Feller).

\begin{theorem}\label{wf} The transition kernel $P(d\pi_{t+1}|\pi_t ,Q_t ,\eta_t)$ has the weak Feller property, as defined in \cite[C.3]{HernandezLermaMCP}. That is, it is weakly continuous in $\cal{P}(\X)\times\cal{Q}\times\cal{H}$ in the sense that for any $g\in C_b(\mathcal{P(\X)})$,   
    $$ \int_{\cal{P}(\X)} g(\pi_{t+1})P(d\pi_{t+1}|\pi_t,Q_t,\eta_t) \in C_b(\mathcal{P(\X)}\times\mathcal{Q}\times\mathcal{H}).$$
\end{theorem}  
\textbf{Proof.}
Fix some $t \ge 0$ and let $(Q^n_t)_{n \ge 0}$ and $(\eta^n_t)_{n \ge 0}$ be two sequences such that $Q^n_t\rightarrow Q_t$ and $\eta^n_t\rightarrow \eta_t$. Since $\cal{Q}$ and $\cal{H}$ are finite, this means that there exist some $N$ and $M$ such that for all $n \geq N$, $Q_t^n = Q_t$ and for all $n\geq M$, $\eta_t^n = \eta_t$. Thus, for a function $g$ on $C_b(\cal{P}(\X))$ we have that for $n\geq\max(N,M)$,
\begin{align*}
    E[g(\pi_{t+1})|\pi_t^n,Q^n_t,\eta^n_t] = E[g(\pi_{t+1})|\pi_t^n,Q_t,\eta_t].
\end{align*}

With this, the proof then closely follows \cite[Lemma 6 and Lemma 11]{YukLinZeroDelay}. Consider now: 
\begin{align*}
    & E[g(\pi_{t+1})|\pi_t,Q_t,\eta_t] \\ 
    &=  \sum_{q_t\in\cal{M}} g(\pi_{t+1}) P(q_t | \pi_t, Q_t, \eta_t) =  \sum_{q_t\in\cal{M}} g(\pi_{t+1})\pi_t(Q_t^{-1}(q_t)) \\
    &=  \sum_{q_t\in\cal{M}} g(F(\pi_t,Q_t,\eta_t,q_t))\pi_t(Q_t^{-1}(q_t)).
\end{align*}
Let $\pi_t^n\rightarrow\pi_t$ weakly. That is, for all continuous and bounded $h$, we have:
$$\int hd\pi_t^n\rightarrow\int hd\pi_t.$$
Note that since $\X$ is finite, this is equivalent to $\pi^n_t \to \pi_t$ in total variation, i.e.,
\begin{equation*}
    \sup_{||h||_\infty \le 1} \left| \int h d\pi^n_t - \int h d\pi_t \right| \to 0,
\end{equation*}
where the supremum is over all measurable functions bounded by $1$. Now, we have that
\begin{align}
& \quad \: \Bigl| \sum_{q_t\in\cal{M}} g(\pi_{t+1})\pi_t(Q_t^{-1}(q_t)) - \sum_{q_t\in\cal{M}} g(\pi^n_{t+1})\pi^n_t(Q_t^{-1}(q_t)) \Bigr| \nonumber \\
& \leq \Bigl| \sum_{q_t\in\cal{M}} g(\pi_{t+1})\pi_t(Q_t^{-1}(q_t)) - g(\pi_{t+1})\pi^n_t(Q_t^{-1}(q_t)) \Bigr| \nonumber \\
& \qquad + \Bigl| \sum_{q_t\in\cal{M}} g(\pi_{t+1})\pi^n_t(Q_t^{-1}(q_t)) - g(\pi^n_{t+1})\pi^n_t(Q_t^{-1}(q_t)) \Bigr|,\label{eq:WF_triangle}
\end{align}
where $\pi^n_{t+1} = F(\pi^n_t, Q_t, \eta_t, q_t)$. By the fact that $\pi_t^n\rightarrow\pi_t$ in total variation, we have that $\pi^n_t(Q^{-1}_t(q_t)) \to \pi_t(Q^{-1}_t(q_t))$, and thus the first term goes to $0$.

For the second term, we have that
\begin{align*}
    \pi^n_{t+1}(\cdot) & = F(\pi_t^n,Q_t,\eta_t,q_t) \\ & = \frac{1}{\pi^n_t(Q_t^{-1}(q_t))} \sum_{x_t\in Q_t^{-1}(q_t)} P(\cdot|x_t,\eta_t(Q_t,q_t))\pi^n_t(x_t),
\end{align*}
which, by the total variation convergence of $\pi^n_t$, converges to
\begin{equation*}
    \frac{1}{\pi_t(Q_t^{-1}(q_t))} \sum_{x_t\in Q_t^{-1}(q_t)} P(\cdot|x_t,\eta_t(Q_t,q_t))\pi_t(x_t) = \pi_{t+1}(\cdot)
\end{equation*}
in total variation for every $q_t$ such that $\pi_t(Q_t^{-1}(q_t)) > 0$. Since $g$ is continuous, this implies that the second term in \eqref{eq:WF_triangle} converges to 0, and the result follows.\qed

Note that a version of this theorem (in the non-controlled case) was proven in \cite[Lemma 11]{YukLinZeroDelay}, but for more general sources. Similar arguments could be applied here (see \cite[Lemmas 3, 6, 11]{YukLinZeroDelay}) but since we are only considering the finite source/action setup here, this proof is more direct.

\subsection{Cost Equivalence and Optimality of Controlled-Predictor-Structured Policies}


By Theorem \ref{Walrand and Varaiya}, for any $\overline{\gamma} \in \Gamma_A$, there exists some $\overline{\gamma}^* \in \Gamma_{C-P}$ such that
$$E^{\overline{\gamma}^*}_{\pi_0} \left[\sum^{N-1}_{k=0} c(x_k,u_k) \right] \le E^{\overline{\gamma}}_{\pi_0} \left[\sum^{N-1}_{k=0} c(x_k,u_k) \right].$$
Now we have that
 \begin{align*}
       &E^{\overline{\gamma}^*}_{\pi_0} \left[\sum^{N-1}_{k=0}c(x_k,u_k)\right] \\&=  E^{\overline{\gamma}^*}_{\pi_0} \left[\sum^{N-1}_{k=0}E^{\overline{\gamma}^*}_{\pi_0}\left[c(x_k,u_k)\mid q_{[0,k-1]}, Q_k, \eta_k \right]\right]\\
      & = E^{\overline{\gamma}^*}_{\pi_0} \left[\sum^{N-1}_{k=0} \sum_{\X\times\mathbb{U}\times\M} \eta_k(u_k|q_k,Q_k)Q_k(q_k|x_k)\pi_k(x_k)c(x_k,u_k) \right] \\
      & = E^{\overline{\gamma}^*}_{\pi_0} \left[\sum^{N-1}_{k=0} \tilde{c}(\pi_k, Q_k, \eta_k) \right],
\end{align*}
where
\begin{align} 
    \tilde{c}(&\pi_k,Q_k,\eta_k) \nonumber
    \\&\coloneqq \sum_{\X\times\mathbb{U}\times\M} \eta_k(u_k|q_k,Q_k)Q_k(q_k|x_k)\pi_k(x_k)c(x_k,u_k)\label{new cost funtion}
    \end{align}
and the second equality follows from the fact that under $\overline{\gamma}^* \in \Gamma_{C-P}$, $Q_k$ and $\eta_k$ are deterministic functions of $\pi_k$ (and hence of $q_{[0,k-1]}$). To be succinct, we have used the notation $Q_k(q_k | x_k)$ and $\eta_k(u_k|q_k,Q_k)$ for the ``measures'' induced by these deterministic functions, but of course since these are all Dirac measures it would be equivalent to rewriting the sum over the relevant sets.

Using the Markov property of $(\pi_k,Q_k,\eta_k)$ from Theorem \ref{Markov property}, we have thus a finite horizon cost criterion for a reformulated problem as
\begin{equation}\label{NewMDPCost}
    J^N(\pi_0, \overline{\gamma}) = E^{\overline{\gamma}}_{\pi_0}\left[\sum^{N-1}_{k=0} \tilde{c}(\pi_k,Q_k,\eta_k)\right]
\end{equation}
for some $\overline{\gamma} \in \Gamma_{C-P}$, where $\tilde{c}$ is the equivalent cost function at each \textit{effective} state, $\pi_t$, and \textit{effective} action, $(Q_t,\eta_t)$. Further, let us define the minimum finite horizon cost criterion as
\begin{equation}\label{mincost}
    J^N(\pi_0) := \inf_{\overline{\gamma}\in \Gamma_{C-P}} E^{\overline{\gamma}}_{\pi_0}\left[\sum^{N-1}_{k=0} \tilde{c}(\pi_k,Q_k,\eta_k)\right].
\end{equation}

This is now a standard Markov Decision Problem (MDP) which then admits an optimal policy (of Markov type) given the weak Feller condition (see Theorem \ref{wf}) under measurable selection conditions for weakly continuous kernels \cite[Chapter 3]{HernandezLermaMCP}. Thus we have proven the following:

\begin{theorem}\label{costeqiv and blackwell}
The minimum cost of the two problems (\ref{objective}) and (\ref{NewMDPCost}) are equivalent. That is,
    \begin{equation*} 
        \inf_{\overline{\gamma}\in\Gamma_A} E^{\overline{\gamma}}_{\pi_0} \left[\sum^{N-1}_{k=0} c(x_k,u_k) \right] = \inf_{\overline{\gamma}\in\Gamma_{C-P}}E^{\overline{\gamma}}_{\pi_0} \left[\sum^{N-1}_{k=0} \tilde{c}(\pi_t,Q_t,\eta_t)\right].
    \end{equation*}
    Accordingly, an optimal coding and control policy in $\Gamma_{C-P}$, when it exists, is also optimal among all admissible coding and control policies.
\end{theorem}

\subsection{Infinite Horizon Discounted Cost Criterion and Existence of Optimal Policies} \label{infinite cost section}

The infinite horizon discounted cost criteria for the reformulated problem is 
    \begin{equation}\label{infinitediscountcost}
J_\beta(\pi_0,\overline{\gamma})=\lim_{N\rightarrow\infty} J_\beta^N(\pi_0,\overline{\gamma})
    \end{equation}
where
\[J_\beta^N(\pi_0,\overline{\gamma})= E^{\overline{\gamma}}\left[\sum_{k=0}^{N-1}\beta^k \tilde{c}(\pi_k,Q_k,\eta_k)\right],\]
and $\beta\in(0,1)$ is the discount factor. The infinite horizon discounted cost problem is to find
 \begin{equation}\label{mininfinitecost}
        J_\beta(\pi_0) = \inf_{\overline{\gamma}\in \Gamma_A} J_\beta(\pi_0,\overline{\gamma}).
    \end{equation}
        We note that a stationary policy is one which maps the Markovian state to actions (consistent with the MDP  formulation) and which does not depend on $t\geq0$. 
    \begin{theorem} \label{best performance class}
        For the system in Figure \ref{loop} and optimization objective (\ref{mininfinitecost}), there exists an optimal policy in $\Gamma_{C-P}$ which is stationary. 
    \end{theorem}
    \textbf{Proof.}

        From Theorem \ref{costeqiv and blackwell}, we know that an optimal policy for \eqref{mincost} will be in $\Gamma_{C-P}$. Using the following inequality, we will now argue that a lower bound for \eqref{mininfinitecost} can be achieved by a minimizing stationary policy for \eqref{mincost}.
    \begin{align}
        J_\beta(\pi_0) = \inf_{\overline{\gamma}\in\Gamma_A}\lim_{N\rightarrow\infty} J_\beta^N(\pi_0,\overline{\gamma})\geq \lim_{N\rightarrow\infty} \inf_{\overline{\gamma}\in\Gamma_A} J^N_\beta(\pi_0,\overline{\gamma})
    \end{align}
    First note that:
    \begin{align*}
         J_\beta^N &(\pi_0,\overline{\gamma}) = \\
         &E_{\pi_0}^{\overline{\gamma}}[\tilde{c}(\pi_0,Q_0,\eta_0)+\beta E[\sum_{k=1}^{N-1}\beta^{k-1}\tilde{c}(\pi_k,Q_k,\eta_k)]|\pi_0,Q_0,\eta_0]
    \end{align*}
    Taking the infimum over $\overline\gamma\in\Gamma_A$ and noting that for finite horizons $\Gamma_{C-P}$ is an optimal class:
    \begin{align*}
        J_\beta^N (\pi_0) &= \inf_{\overline\gamma\in\Gamma_A} E_{\pi_0}^{\overline{\gamma}}[\tilde{c}(\pi_0,Q_0,\eta_0)+\beta E[J_\beta^{N-1}(\pi_1)|\pi_0,Q_0,\eta_0] \\
        & = \inf_{\overline\gamma\in \Gamma_{C-P}} E_{\pi_0}^{\overline{\gamma}}[\tilde{c}(\pi_0,Q_0,\eta_0)+\beta E[J_\beta^{N-1}(\pi_1)|\pi_0,Q_0,\eta_0]
    \end{align*}
    As $N$ increases, $J_\beta^N (\pi_0)$ is also increasing monotonically and by a contraction argument \cite{HernandezLermaMCP} we have that $\lim_{N\rightarrow\infty} J_\beta^N (\pi_0) = J^\infty_\beta(\pi_0)$ exists,
    \begin{align*}
        J^\infty_\beta(\pi_0) \coloneq\lim_{N\rightarrow\infty} \min_{Q_0,\eta_0} (\tilde{c}(\pi_0,Q_0,\eta_0)+\beta E[J^{N-1}_\beta(\pi_1)|\pi_0,Q_0,\eta_0])
    \end{align*}
which thus serves as a lower bound to the optimal cost. 

    The cost function $c(x,u)$ is bounded and since $\tilde{c}(\pi_t, Q_t, \eta_t)$ is continuous in $\pi_t$, we have that $\tilde{c}(\pi_t,Q_t,\eta_t)$ is continuous and bounded on $\mathcal{P}(\X)\times\mathcal{Q}\times\mathcal{H}$. Furthermore, $\cal{Q} \times \cal{H}$ is compact, and by Theorem \ref{wf} we have that $P(d\pi_{t+1} | \pi_t, Q_t, \eta_t)$ is weak Feller. Thus (see \cite[Chapter 8.5]{hernandezlasserre1999further}) there exists a stationary $\overline{\gamma} \in \Gamma_{C-P}$ which is optimal for \eqref{infinitediscountcost} (that is, it satisfies \eqref{mininfinitecost}). For a concise discussion, see also \cite[Lemma 5.5.4]{yuksel2020control}.   
    \qed

\subsection{Interpretation of the Structural Results in Special Cases}
The problem setup being considered is one which is dynamically controlled and in which the measurements are obtained via coding. This setup is a generalization of the following two important settings.
\begin{enumerate}
    \item {\bf Zero-Delay Coding of Control-Free Markov Sources.} Consider the model (\ref{system}) where control is absent. Dynamics for this system are described by $x_{t+1} = f(x_t,w_t)$, for all $t \geq 0$. In this setup, the set of admissible encoder policies is 
    $$\gamma_t^e:I_t^e\to \cal{M},$$
    where $I^e_t = (x_{[0,t]}, q_{[0,t-1]})$. Instead of a controller we would have a decoder and a decoder policy, called $\gamma^d$. The decoder policy would map the channel output to a reconstruction of the source $\hat{x}_t$ defined on $\X$. That is, 
    \begin{align*}
        \gamma^d_t: {\cal{M}}^{t+1} \to \X.
    \end{align*}
    An optimization objective for this problem is to minimize the difference between the state and the reconstructed state. Let $\gamma^e = \{\gamma^e_t\}_{t\geq0}$ and $\gamma^d = \{\gamma^d_t\}_{t\geq0}$. The finite horizon cost is:
    \begin{equation*}
        J^N (\pi_0,\gamma^e, \gamma^d) = E_{\pi_0}^{\gamma^e,\gamma^d}\Bigg[\sum^N_{t=0} d(x_t,\hat{x}_t)\Bigg]
    \end{equation*}
    where $d:\X\times\hat{\X}\mapsto\mathbb{R}_+$ is a distance metric. The reconstructed state is not fed back into the plant and does not impact the next state value, $x_{t+1}$. We can view this problem as a special case of our framework by letting $u_t = \hat{x}_t$ and $x_{t+1} = f(x_t, u_t, w_t) \equiv f(x_t, w_t)$ for all $u_t$. Our ``cost'' would then become $c(x,u) = d(x, \hat{x})$.

    By Theorem \ref{Walrand and Varaiya}, we know that the optimal encoding-decoding policy is of the form $\pi_t \to (Q_t, \eta_t)$. In the non-controlled case, we can explicitly identify the optimal decoder for a given encoder by:
    $$\gamma^d_t(q_{[0,t]}) = \min_{\hat{x}} E_{\pi_0}^{\gamma^e}\left[ d(x_t, \hat{x}) | q_{[0,t]} \right],$$
    which indeed has the form $\gamma^{c}_t:\pi_t\mapsto \Bigl( \eta_t:(Q_t ,q_t)\mapsto u_t \Bigr)$, as in \eqref{pi_to_eta}. Similarly, it was previously shown in \cite{WalrandVaraiya} that the optimal encoder has the form
    $$\gamma^{e}_t:\pi_t\mapsto (Q_t: x_t\mapsto q_t),$$
    as in \eqref{pi_to_Q}.
    
    Since the decoder/controller can be explicitly identified in this case, and since the ``control'' doesn't affect the source evolution, this can instead be studied as a Markov chain $(\pi_t,Q_t)$ to only search for the optimal encoder $\gamma^e$. This leads to a simpler setup as studied in \cite{cregg2023reinforcement,wood2016optimal}.

    \item {\bf Partially Observed Markov Decision Problems (POMDPs).}

 Instead of using an encoding policy, $\gamma^e_t$, to find a quantization of $x_t$, suppose that the quantization action is fixed. 
    Such a fixed quantization induces a partially observed Markov decision problem (POMDP). The coding policy is a constant-valued  function:
    \begin{equation*}
        \gamma^e_t:\pi_t\mapsto \tilde{Q} .
    \end{equation*}
    This is equivalent to a system with dynamics:
    \begin{align*}
        x_{t+1} &= f(x_t,u_t,w_t)\\
        q_t &= \tilde{Q}(x_t)
    \end{align*}
    The controller will only have access to channel output $q_t$ making it analogous to an observation of $x_t$ at time t. When the control problem is structured in this form, the problem is reduced to that of a POMDP with observation $q_t$ of $x_t$. Structural results for POMDPs are well established  \cite{Yus76,Rhe74,Blackwell2} and it is well-known that an optimal control will be using the filter $\overline{\pi}_t$ (see (\ref{filter_2})) as a sufficient statistic; here we present an alternative interpretation which is consistent with the real-time encoding and control framework: the control policy will be structured as
    \begin{equation*}
    \gamma^c_t:\pi_t \mapsto (\eta_t: q_t\mapsto u_t )
    \end{equation*}
    or 
    \begin{equation*}
    \gamma^c_t:(\pi_t,q_t) \mapsto u_t.
    \end{equation*}
    Note that $\overline{\pi}_t$ is determined by $(\pi_t,q_t)$ in this setup, so indeed POMDPs and the optimality of filter policies are a special case of our problem.
    \end{enumerate}

\section{Finite State Approximation via Predictor Quantization, its Near Optimality, and Reinforcement Learning} \label{section:Predictor_Quantization} 
Due to the uncountable state space of probability measures on $\cal{P}(\X$), computation of value functions is challenging. To address this problem, we will quantize this uncountable space to arrive at an approximate finite space model whose solution will be near-optimal for the original, building on the analysis introduced in \cite[Chapter 4]{SaLiYuSpringer} for weakly continuous MDPs. Note that this quantization is different than the {\it operational} quantization introduced in Section \ref{problemsetup} where a finite space $\X$ is quantized into a smaller alphabet due to physical constraints in networked control.
\subsection{Finite MDP Approximation and its Near Optimality}
Quantization (or discretization) of an MDP on a possibly uncountable space $\Z$ will be done as follows. We will define ${\Z_n}=\{\hat{z}_{n,1},...,\hat{z}_{n,m_n}\}$ as a finite set that approximates $\Z$ and a measurable mapping $\rho:\Z \mapsto \Z_n$. If we define a disjoint partition $\{B_{n,i}\}_{i=1}^{m_n}$ of $\Z$, and pick a representative $\hat{z}_{n,i}$ in each $B_{n,i}$ such that
$$\rho(z) = \hat{z}_{n,i} \quad \forall z\in B_{n,i}$$ 
then $\Z_{n}$ can be considered a finite set of representative states for bounded sections of uncountable $\Z$ \cite{SaLiYuSpringer}. More specifically, we will define $\rho(z)=\hat{z}$ as the nearest neighbour map:
\begin{equation*}
    \hat{z} = \argmin_{z'\in \Z_n} d_z(z,z')
\end{equation*}
where $d_z(\cdot,\cdot)$ is a metric on $\Z$. 

Building on \cite[Chapter 4]{SaLiYuSpringer}, we construct an approximate finite model and to establish near optimality of this approximate MDP we will follow \cite[Section 2.3]{KSYContQLearning}. First, we will select a measure $\kappa \in \mathcal{P}(\Z)$, for which $\kappa(B_{n,i})\geq 0$ for all $B_{n,i}$. Using this measure, we can define a stage-wise cost and transition kernel. We will denote the cost and kernel from the original MDP as $c(z,u)$ and $P(dz'|z,u)$, respectively. For $\hat{z}_{n,i}, \hat{z}_{n,j} \in \Z_n$ and $u \in \mathbb{U}$ we define:
\begin{equation} \label{quantized cost}
    c_n(\hat{z}_{n,i}, u) \coloneqq \int_{B_{n,i}} \frac{\kappa(dz)}{\kappa(B_{n,i})} c(z,u)
\end{equation}
\begin{align}
    P_n&(\hat{z}_{n,j}|\hat{z}_{n,i}, u) \label{quantized kernel}
    \coloneqq \int_{z\in B_{n,i}}\int_{z'\in B_{n,j}} P(dz' \mid z, u)\frac{\kappa(dz)}{\kappa(B_{n,i})}
\end{align}
We can now define our approximate MDP:
\begin{definition}
For an MDP with state space $\Z$, action space $\mathbb{U}$, transition kernel $P(dz'|z,u)$, and cost function $c(z,u)$, which we denote by $\text{MDP} = (\Z, \mathbb{U}, P(dz'|z,u), c(z,u))$ we define the \textit{Finite State Approximate MDP} as $\text{MDP}_n \coloneqq (\Z_n, \mathbb{U}, P_n(\hat{z}_{n,j} | \hat{z}_{n,i},u), c_n(\hat{z}_{n,i},u))$.
\end{definition}
The optimal value function for $\text{MDP}_n$ is the solution $$\hat{J}_\beta:\Z_n \mapsto \R$$ to the Discounted Cost Optimality Equation (DCOE):
$$\hat{J}_\beta(\hat{z}_{n,i}) = \min_{u\in\mathbb{U}} (c_n(\hat{z}_{n,i},u)+\beta \sum_{\hat{z}_{n,j}\in\Z_n} \hat{J}_\beta(\hat{z}_{n,j})P_n(\hat{z}_{n,j}\mid \hat{z}_{n,i},u))$$
for $\hat{z}_{n,i} \in \Z_n$. Note that we can also extend this function over $\Z$ as follows: if $\hat{z}_{n,i}\in B_{n,i}$, then for all $z\in B_{n,i}$ we have $\hat{J}_\beta (z) \coloneqq \hat{J}_\beta (\hat{z}_{n,i})$ \cite{SaLiYuSpringer}. Equivalently, we can extend a policy $\hat{\gamma}_n$ for the $\text{MDP}_n$ to a policy $\tilde{\gamma}$ on the original MDP as follows. 
\begin{equation}
    \tilde\gamma(z) = \hat{\gamma}_n(\hat{z}) \label{extension over MDP}
\end{equation}
wherever $\hat{z} = \rho(z)$. 
We now state an assumption on the original MDP.
\begin{assumption}\cite[Section 2.3]{KSYContQLearning}\cite[Chapter 4]{SaLiYuSpringer}
\label{quantization assumption}
    \begin{enumerate}[(i)]
        \item The cost function $c(z,u)$ is continuous and bounded. 
        \item The transition kernel $P(dz'|z,u)$ is weakly continuous.
        \item Spaces $\Z$ and $\mathbb{U}$ are compact. 
    \end{enumerate}
\end{assumption}
Using the compactness of $\Z$ from this assumption, we have that
\begin{equation*}
    \lim_{n\rightarrow\infty} \max_{z\in\Z} \min_{i = 1,...,m_n} d_z(z,\hat{z}_{n,i}) = 0
\end{equation*}
which motivates the following theorem.
\begin{theorem} \cite[Theorem 4.3]{SaLiYuSpringer}\\ \label{quantization performance}
    Let $\tilde\gamma^*$ be a policy extended over the MDP as in \eqref{extension over MDP} from an optimal $\hat{\gamma}_n^*$ on $\text{MDP}_n$. Then for all $z_0\in \Z$ we have $$\lim_{n\rightarrow\infty} |J_\beta(z_0,\tilde\gamma^*) - J_\beta(z_0)| = 0.$$
\end{theorem}
That is, the optimal policy for $\text{MDP}_n$, when appropriately extended over $\Z$, becomes near-optimal for the original MDP.


Consider now the Controlled-Predictor MDP introduced in Section \ref{Predictor Structured Markov Problem}; that is, $\text{MDP} = (\mathcal{P}(\X), \mathcal{Q} \times \mathcal{H}, P(d\pi' | \pi, Q, \eta), \tilde{c}(\pi, Q, \eta))$. We first define the quantization of the belief-space $\cal{P}(\X)$ to a finite $\mathcal{P}_n(\X)$, where $n$ represents the resolution of the quantization. 
Let $\mathcal{P}_n(\X)$ be the following finite set:
\begin{align*}
    \{\hat{\pi}\in \mathcal{P}(\X)\colon \hat{\pi}=\left[\frac{k_1}{n},...,\frac{k_{|\X|}}{n}\right],k_i \in \{0,...,n\},i = 1,...,|\X|\}
\end{align*}
Using the nearest neighbour quantization specified above, we define $\rho$:
\begin{equation*}
    \rho(\pi) \coloneq \argmin_{\hat{\pi}\in \mathcal{P}_n(\X)} d(\pi, \hat{\pi}) 
\end{equation*}
where $d(\cdot,\cdot)$ is a metric on $\mathcal{P}(\X)$ (say the $L^1$ distance, since $\X$ is finite) \cite{SaLiYuSpringer}.
The partition of $\mathcal{P}(\X)$ that is induced by $\rho$ is $$\{B_{n,i}\} = \{\pi \in \mathcal{P}(\X): \rho(\pi) = \hat{\pi}_{n,i}\}.$$ 
We define $c_n(\hat{\pi}_{n,i},Q,\eta)$ and $P_n(\hat{\pi}_{n,j} | \hat{\pi}_{n,i}, Q, \eta)$ analogously to \eqref{quantized cost} and \eqref{quantized kernel}, and define $\text{MDP}_n \coloneqq (\mathcal{P}_n(\X), \mathcal{Q} \times \mathcal{H}, P_n(\hat{\pi}_{n,j} | \hat{\pi}_{n,i}, Q, \eta), c_n(\hat{\pi}_{n,i}, Q, \eta))$.

We know, from Theorem \ref{wf}, that the transition kernel $P(d\pi' | \pi, Q, \eta)$ is weak Feller continuous. We have also discussed the boundedness and continuity of $\tilde{c}(\pi,Q,\eta)$ in the proof of Theorem \ref{best performance class}. Spaces $\cal{Q}$ and $\cal{H}$ are compact (finite). Since $\X$ is finite, $\cal{P}(\X)$ is compact, so we have satisfied Assumption \ref{quantization assumption}. Therefore, from Theorem \ref{quantization performance}, we have the following result concerning the optimality of the solution to $\text{MDP}_n$.
\begin{theorem}\label{special case quantization performance}
If the optimal policy $\hat{\gamma}^*_n$ for $\text{MDP}_n$ is extended to $\tilde\gamma_n^*$ on the original controlled-predictor MDP, 
we have for all $\pi_0\in \cal{P}(\X)$:
\begin{equation}
    \lim_{n\to\infty} \left|J_\beta(\pi_0,\tilde\gamma_n^*) - J_\beta(\pi_0)\right|  = 0
\end{equation} 
\end{theorem}

\subsection{Quantized Q-Learning}
We will now consider a learning algorithm based on Watkins and Dayan's Q-Learning algorithm \cite{Watkins} and its extension to Quantized Q-Learning  \cite{YukselBasarBook24}. Noting the above near optimality result, we will show that we can run the Q-learning algorithm on the approximate $\text{MDP}_n$ to achieve a near optimal result for the original. We will first comment on the convergence of Q-Learning for non-Markovian environments. The following is a result from \cite{karayukselNonMarkovian}. Let us have the following sequences:
\begin{enumerate}
    \item $\{y_t\}_{t\geq 0}$ defined on finite space $\mathbb{Y}$.
    \item $\{u_t\}_{t\geq 0}$ defined on finite space $\mathbb{U}$.
    \item $\{C_t\}_{t\geq 0}$ defined on $\mathbb{R}$.
    \item $\{\alpha_t\}_{t\geq 0}$, where $\alpha_t(y,u):\mathbb{Y}\times\mathbb{U} \rightarrow \mathbb{R}$ is the \emph{learning rate} at time $t$.
    \item $\{\mathbf{Q}_t\}_{t\geq 0}$, where $\mathbf{Q}_t:\mathbb{Y}\times\mathbb{U} \rightarrow \mathbb{R}$ and $\mathbf{Q}_0 \equiv 0$, are the \emph{Q-factors} at time $t$.
\end{enumerate}
We define the following Q-learning iterations for every state-action pair $(y,u)\in\mathbb{Y}\times\mathbb{U}$:
\begin{align}
    \mathbf{Q}_{t+1}(y,u) =& (1-\alpha_t(y,u))\mathbf{Q}_t(y,u) \nonumber\\
    &+ \alpha_t(y,u)(c_t+\beta\min_{u\in\mathbb{U}}\mathbf{Q}(y_{t+1},u))
\end{align}
Under the following assumption, we have a convergence theorem. 
\begin{assumption}\label{Qconvergence assumption}\cite{karayukselNonMarkovian}\\
    \begin{enumerate}[(i)]
        \item $\alpha_t(y,u) = 
            \begin{cases} 
                \frac{1}{1+\sum^t_{k=0} \mathds{1}_{\{y_k=y,u_k=u\}}} & \text{if } (y_t,u_t) = (y,u) \\
                0 & \text{otherwise} \\
            \end{cases}$
            \\and $\sum_{t=0}^\infty \alpha_t(y,u) = \infty$ almost surely.
        \item For some function $c^*:\mathbb{Y}\times\mathbb{U}\mapsto\R$, we have that almost surely
        $$\frac{\sum^t_{k=0} c_k \mathds{1}_{\{y_k=y,u_k=u\}}}{\sum^t_{k=0}  \mathbb{1}_{\{y_k=y,u_k=u\}}} \rightarrow c^*(y,u).$$
        \item For any $g : \mathbb{Y} \to \R$, there is a measure $P^*$ such that almost surely
        $$\frac{\sum^t_{k=0} g(y_{k+1}) \mathbb{1}_{\{y_k=y,u_k=u\}}}{\sum^t_{k=0}\mathds{1}_{\{y_k=y,u_k=u\}}} \rightarrow \int g(y') P^*(dy'|y,u).$$
    \end{enumerate} 
    \end{assumption}
    Parts (ii) and (iii) follow from assuming that each $(y,u)\in(\mathbb{Y}, \mathbb{U})$ is visited infinitely often. This is a standard reinforcement learning algorithm ergodic-type condition and is especially practically attainable and nonrestrictive and when $\mathbb{Y}$ and $\mathbb{U}$ are finite.
\begin{theorem} \cite[Theorem 2.1]{karayukselNonMarkovian}. \label{Q-learning convergence theorem}
\begin{enumerate}
    \item $\mathbf{Q}_t(y,u) \rightarrow \mathbf{Q}^*(y,u)$ almost surely for each $(y,u)\in\mathbb{Y}\times\mathbb{U}$.
    \item $\mathbf{Q}^*(y,u)$ is the solution to $$\mathbf{Q}^*(y,u) = c^*(y,u) +\beta\sum_{y'\in\mathbb{Y}} \min_{u'\in\mathbb{U}} \mathbf{Q}^*(y',u') P^*(y'|y,u) $$
    \item An optimal policy for $\text{MDP}_n = (\mathbb{Y}, \mathbb{U}, P^*, c^*)$ is given by: $\hat{\gamma}_n^*(y)\coloneq \argmin_{u\in\mathbb{U}} \mathbf{Q}^*(y,u)$.
    \end{enumerate}
\end{theorem}
In order to apply the above result to the MDP in the previous sections, we require the following assumption:
\begin{assumption}\label{assumption:PHR}
The transition kernel $P(x'|x,u)$ is such that the induced Markov chain on the state space for a corresponding fully observable Markov Decision Process under any stationary policy is irreducible and aperiodic.
\end{assumption}
    Under Assumption \ref{assumption:PHR}, if $(Q_t, \eta_t)$ is chosen independently and uniformly from $\mathcal{Q} \times \mathcal{H}$ at each $t \ge 0$, the resulting Markov process $\{x_t\}_{t \ge 0}$ is positive Harris recurrent (see e.g., \cite[Definition 3.2.4]{yuksel2020control}), and thus admits a unique invariant probability measure $\zeta$.

We apply the above algorithm to our (approximate) controlled-predictor MDP. Accordingly, let $y_t=\hat{\pi}_t, u_t = (Q_t,\eta_t),$ and $C_t=\tilde{c}(\pi_t,Q_t,\eta_t)$.
To show that this algorithm converges, we will need to satisfy Assumption \ref{Qconvergence assumption} for the sequence $(\hat{\pi}_t,(Q_t,\eta_t),\tilde{c}(\pi_t, Q_t, \eta_t), \alpha_t) = (y_t,u_t,C_t,\alpha_t)$.

We will now present some results to validate these assumptions. We follow the strategy in \cite[Section IV]{creggZeroDelayNoiseless} to arrive at the unique ergodicity of the predictor process. Indeed, through a uniformly and independently chosen $Q_t$ and $\eta_t$, the system would be equivalent to a control-free hidden Markov model. Therefore by \cite[Lemma 4]{creggZeroDelayNoiseless} we have that the predictor process is stable in total variation and by \cite[Theorem 4]{creggZeroDelayNoiseless}, we have unique ergodicity. Due to the specific nature of the problem, we also have the following additional properties:

First, we have that there exists some $Q\in\cal{Q}$, $x\in\X$, $q\in \mathcal{M}$ such that $Q$ quantizes $x$ without any loss to $q$ (that is, $Q^{-1}(q) = \{x\}$). We also have, by the finiteness of $\cal{Q}$ and $\cal{H}$ and the positive Harris recurrence of $\{x_t\}_{t \ge 0}$ (under Assumption \ref{assumption:PHR}), that for some $t\geq0$ and all $\eta \in \cal{H}$, we have $(x_t,Q_t,\eta_t) = (x,Q,\eta)$ almost surely \cite{creggZeroDelayNoiseless}. This implies, through the update equation \eqref{predictor update}, that
\begin{align}
    \pi_{t+1}(x') &= \frac{1}{\pi_t(Q^{-1}(q))} \sum_{x\in Q^{-1}(q)} P(x'|x,\eta(Q,q))\pi_t(x) \nonumber \\
   &= \frac{P(x'|x,\eta(Q,q))\pi_t(x)}{ \pi_t(x)} \nonumber \\
   &= P(x'|x,\eta(Q,q))\label{recurrent_pi}
\end{align}
 
Therefore, predictors of the form in \eqref{recurrent_pi} act as ``atoms'' with an almost surely finite return time, which implies that there is at most one invariant measure for $\{\pi_t\}_{t \ge 0}$ under the uniform selection of $(Q_t, \eta_t)$. Existence of an invariant measure is guaranteed by the compactness of $\mathcal{P}(\X)$ and the weak Feller property of $P(d\pi_{t+1} | \pi_t, Q_t, \eta_t)$ (see \cite[Theorem 7.2.3]{Hernandez-Lerma2003}). 

Thus, our process $\{\pi_t\}_{t \ge 0}$ admits a unique invariant measure, which we denote by $\lambda$. This leads to the following:
\begin{lemma} \label{ergodic}
        Assume that at each $k \ge 0$, $(Q_k, \eta_k)$ is chosen uniformly from $\mathcal{Q} \times \mathcal{H}$, and let Assumption \ref{assumption:PHR} hold. Then for any measurable and bounded function $g:\mathcal{P}(\X)\mapsto\R$ and any $\pi_0 \in P(\X)$,
        \begin{equation*}
            \frac{1}{t}\sum^{t-1}_{k=0} g(\pi_k) \rightarrow\int g d\lambda
        \end{equation*}
\end{lemma}
\textbf{Proof.}
    By the pathwise ergodic theorem \cite[Corollary 2.5.2]{Hernandez-Lerma2003}, and the fact that for any $\pi_0$, the hitting time to $\pi_t$ of the form in \eqref{recurrent_pi} is almost surely finite, the result follows. \qed
    
We will now define:
$$B_n^\lambda \coloneqq \{B\in \{B_{n,i}\}:\lambda(B)>0\}$$
and
$$\mathcal{P}_n^\lambda(\X) \coloneqq \{\hat{\pi}\in \mathcal{P}_n(\X): \rho^{-1}(\hat\pi)\in B_n^\lambda\}.$$
We then have the following:
\begin{lemma} For any $\pi_0$ and $(\hat\pi,Q,\eta)\in \mathcal{P}^\lambda_n(\X)\times\cal{Q}\times{H}$, under an independent and uniformly distributed $(Q_k, \eta_k)_{k \ge 0}$, we have that almost surely:
\begin{enumerate}[(i)]
    \item $(\hat\pi_t,Q_t,\eta_t) = (\hat\pi,Q,\eta)$ infinitely often and $\sum_{t\geq0} \alpha_t(\hat\pi,Q,\eta)=\infty$.
    \item
        $$\frac{\sum^t_{k=0} \tilde{c}(\pi_k, Q_k, \eta_k) \mathds{1}_{\{\hat\pi_k=\hat\pi,Q_k=Q,\eta_k=\eta\}}}{\sum^t_{k=0}  \mathds{1}_{\{\hat\pi_k=\hat\pi,Q_k=Q,\eta_k=\eta\}}} \rightarrow c_n(\hat{\pi},Q,\eta),$$
        where 
        $$c_n(\hat{\pi},Q,\eta) \coloneqq \int_{B} \frac{\lambda(d\pi)}{\lambda(B)}\tilde{c}(\pi, Q, \eta)$$
        and $B$ is the bin of $\hat{\pi}$.
    \item If $P_n(d\hat{\pi}'|\hat{\pi},Q,\eta)$ is defined as in \eqref{quantized kernel} using $\lambda$ instead of measure $\kappa$, we have
    \begin{multline*}
        \frac{\sum^t_{k=0} g(\hat\pi_{k+1}) \mathds{1}_{\{\hat\pi_k=\hat\pi,Q_k=Q,\eta_k=\eta\}}}{\sum^t_{k=0} \mathds{1}_{\{\hat\pi_k=\hat\pi,Q_k=Q,\eta_k=\eta\}}} \\
        \to \int g(\hat\pi') P_n(d\hat\pi'|\hat\pi,Q,\eta).
    \end{multline*}
\end{enumerate}
\end{lemma}
\textbf{Proof.}
    Using Lemma \ref{ergodic} with $(\pi_t,Q_t,\eta_t)$, and noting that $(Q_k, \eta_k)$ are chosen uniformly, we have
    $$\frac{1}{t} \sum^{t-1}_{k=0} g(\pi_k,Q_k,\eta_k) \rightarrow \int g(\pi,Q,\eta) \frac{\lambda}{\abs{\cal{Q}\times\cal{H}}}(d\pi)$$ as $t \rightarrow\infty$, almost surely. 
    Defining
    $$g_1(\pi_k, Q_k, \eta_k) \coloneqq \tilde{c}(\pi_k, Q_k, \eta_k)\mathds{1}_{\{\hat{\pi}_k = \hat{\pi}, Q_k = Q, \eta_k = \eta\}}$$
    and
    $$g_2(\pi_k, Q_k, \eta_k) \coloneqq \mathds{1}_{\{\hat{\pi}_k = \hat{\pi}, Q_k = Q, \eta_k = \eta\}}$$
    we then have
\begin{align*}
    \frac{\sum_{k=0}^t g_1(\pi_k, Q_k, \eta_k)}{\sum_{k=0}^t g_2(\pi_k, Q_k, \eta_k)} \to c_n(\hat{\pi},Q,\eta),
\end{align*}
which satisfies part (ii) of Assumptions \ref{Qconvergence assumption}. Similarly, by letting $g_1, g_2$ be the relevant functions, we also obtain parts (i) and (iii) (see also \cite[Theorem 5.3]{cregg2024reinforcementlearningoptimaltransmission}).\qed

We can now state the following as a result of Theorem \ref{Q-learning convergence theorem} and as an extension of \cite[Theorem 1]{creggZeroDelayNoiseless}.
\begin{theorem} \label{theorem:quantizedQlearning_nearoptimal}
For each $(\hat\pi,Q,\eta) \in \mathcal{P}^\lambda_n(\X)\times\cal{Q}\times{H}$,  $\mathbf{Q}_t (\hat\pi,Q,\eta)$ will converge to a limit almost surely. The limit will satisfy: 
    \begin{align*}
        \mathbf{Q}^*(&\hat{\pi},Q,\eta) = c_n(\hat{\pi},Q,\eta)\\&+\beta\sum_{\hat{\pi}'\in \mathcal{P}_n^\lambda(\X)} P_n(\hat{\pi}'\mid \hat{\pi},Q,\eta) \min_{(Q',\eta')\in\cal{Q}\times\cal{H}} \mathbf{Q}^*(\hat{\pi}',Q',\eta')
    \end{align*}
This coincides with the DCOE for $\text{MDP}_n = (\mathcal{P}^\lambda_n(\X), \mathcal{Q}, P_n, c_n)$, and therefore the policy $\hat{\gamma}^*_n(\hat{\pi}) = \argmin_{Q,\eta}\mathbf{Q}^*(\hat{\pi}, Q, \eta)$ is optimal for $\text{MDP}_n$. Thus if we define a policy $\tilde{\gamma}_n^*(\pi) \coloneq \hat{\gamma}^*_n(\hat\pi)$ where $\hat\pi = \rho(\pi)$, we have by Theorem \ref{special case quantization performance} that:
$$ \lim_{n\to\infty}|J_\beta(\pi_0,\tilde{\gamma}_n^*) - J_\beta(\pi_0)|=0$$ for any $\pi_0$ such that $\lambda(\pi_0) > 0$.
\end{theorem}
\emph{Remark:} One such initialization where $\lambda(\pi_0) > 0$ would be those recurrent predictors identified in \eqref{recurrent_pi}, or any $\pi$ which is reachable from one such $\pi_0$. Since $\pi_{t}$ has only finitely many possible values given $\pi_0$, each of these elements would also have positive measure under $\lambda$.

\section{Sliding Finite Window Structured Policies, Their Near Optimality and Reinforcement Learning}\label{section:Finite_Window}
Instead of the formulation of the networked problem as its equivalent predictor-structured reduction, we will now consider another representation. We will use, as the state, a set of information which includes one previous predictor term, and memory of a finite window of previous values of $Q_t,q_t,$ and $\eta_t$.

To introduce this structure, we will first define a window size $N$ and our new state sequence $\{\kappa_t\}_{t\geq0}$, which is defined on $\W = \mathcal{P}(\X)\times\M^N\times\mathcal{Q}^N\times\mathcal{H}^N$, by
\begin{equation} \label{W}
    \kappa_t = \{\pi_{t-N},I^N_t\},
\end{equation}
where
\begin{equation}
    I_t^N = \{q_{[t-N,t-1]},Q_{[t-N,t-1]},\eta_{[t-N,t-1]}\}.
\end{equation}
$\pi_t$ can be found from $\kappa_t$ by starting from $\pi_{t-N}$ and using the update equation \eqref{predictor update} $N$ times. This map will be denoted as $\psi:\W\mapsto \mathcal{P}(\X)$, so that $\pi_t = \psi(\kappa_t)$.

When at $t<N$, there will not be full window information. For simpler notation, at $t=0$, we will consider $\pi_0$ as found by updating $\pi_{-N}$ $N$ times. 
\begin{definition}
    We will say a joint coder controller policy has \textit{Controlled Finite Sliding Window Structure} if, at time $t$, it uses only $\kappa_t$ to select $Q_t$ and $\eta_t$. We will denote this set of policies as $\Gamma_{C-FW}$.
\end{definition}
As in section \ref{Predictor Structured Markov Problem}, this motivates the definition of a new MDP.
\begin{theorem}
   $(\kappa_t,(Q_t,\eta_t))$ is a controlled Markov chain on $\W\times\cal{Q}\times\cal{H}$.
\end{theorem}
\textbf{Proof.}
\begin{align*}
    & P(\kappa_{t+1} \in \cdot |\kappa_s,Q_s,\eta_s, s \le t) \\
    & = P(\pi_{t-N+1}, q_{[t-N+1, t]}, Q_{[t-N+1, t]}, \eta_{[t-N+1,t]} \in \cdot \\ & \hspace{3cm}| \pi_{[0,t-N]}, q_{[0,t-1]}, Q_{[0,t]}, \eta_{[0,t]}) \\
    & = P(\pi_{t-N+1}, q_{[t-N+1, t]}, Q_{[t-N+1, t]}, \eta_{[t-N+1,t]} \in \cdot \\ & \hspace{3cm}| \pi_{t-N}, q_{[t-N,t-1]}, Q_{[t-N,t]}, \eta_{[t-N,t]}) \\
    & = P(\pi_{t-N+1},I_{t+1} \in \cdot |\pi_{t-N},I_t,Q_t,\eta_t) \\
    & = P(\kappa_{t+1} \in \cdot|\kappa_t,Q_t,\eta_t),
\end{align*}
where the second equality follows from the fact that $\pi_{t-N+1}, q_{[t-N+1,t-1]}, Q_{[t-N+1,t]}$ and $\eta_{[t-N+1,t]}$ are deterministic given $\pi_{t-N}, q_{[t-N,t-1]}, Q_{[t-N,t]},$ and $\eta_{[t-N,t]}$, and also that $Q_t$, $q_t$ depends on the past terms only through $\pi_t$ (which is a function of $\pi_{t-N}, I_t^N$).
\qed

With an abuse of notation, define the following cost function:
\begin{align*}
    c(\kappa&,Q,\eta)
    = \sum_{(x,u,q)\in\X\times\mathbb{U}\times\M} \eta_t(u | Q,q)Q(q|x)\psi(\kappa)(x)c(x,u).
\end{align*}
Thus we obtain $\text{MDP} = (\W, \mathcal{Q} \times \mathcal{H}, P(d\kappa' | \kappa, Q, \eta), c(\kappa,Q,\eta))$. The objective is minimizing the following over all $\overline\gamma\in\Gamma_A$:
\begin{equation}\label{approximate J}
    J_\beta(\kappa_0,\overline\gamma) = E^{\overline\gamma}_{\kappa_0}\Big[\sum^{\infty}_{k=0} \beta^k c(\kappa_k,Q_k,\eta_k) \Big].
\end{equation}
We can extend the results of Theorem \ref{costeqiv and blackwell} and Section \ref{infinite cost section} to argue that $$\inf_{\overline\gamma\in\Gamma_A} J_\beta(x_0,\overline\gamma) = \inf_{\overline\gamma'\in\Gamma_A} J_\beta(\kappa_0,\overline\gamma')$$
and also that an optimal policy for $J_\beta(\kappa_0, \overline{\gamma})$ can be found in $\Gamma_{C-FW}$. That is,
\begin{equation*}
    J_\beta(\kappa_0) = \inf_{\overline\gamma\in\Gamma_{C-FW}}J_\beta(\kappa_0,\overline\gamma) .
\end{equation*}
\subsection{Near Optimality of Finite Sliding Window Approximation}
We will fix the prior $\pi_{t-N}$ to be some (fixed) $\mu$. This approach follows that of \cite[Sections 3.5]{cregg2024reinforcementlearningoptimaltransmission} (see also \cite{kara2020near}). Our state will now be $\hat{\kappa}_t=(\mu,I^N_t)$ defined on $\hat\W=\{\mu\}\times \cal{M}^N \times\cal{Q}^N\times\cal{H}^N$. In this section, we will consider the performance impact of using this approximation. Recalling the definition of $\pi_t$ in \eqref{predictor}, we can define $\pi_t$ starting from the correct prior $\pi_{t-N}$ as:
\begin{equation}\label{finite window predictor}
\pi_t(\cdot) = P^{\overline{\gamma}}_{\pi_{t-N}} (x_t=\cdot |q_{[t-N,t-1]},Q_{[t-N,t-1]}, \eta_{[t-N,t-1]})
\end{equation}
and our approximate predictor that uses $\mu$ as the incorrect prior as:
\begin{equation}\label{finite window approximate predictor}
\hat\pi_t(\cdot) = P^{\overline{\gamma}}_\mu (x_t=\cdot |q_{[t-N,t-1]},Q_{[t-N,t-1]},\eta_{[t-N,t-1]}).
\end{equation}
We will define the following transition kernel for the finite
sliding window approximation. Note that we now make the dependence on $N$ explicit as it affects the quality of the approximation (as will be shown shortly).
\begin{align}\label{finite window approximate kernel}
    P_N(\hat{\kappa}_{t+1}|\hat{\kappa}_t,Q_t,\eta_t) &= P_N(\mu,I_{t+1}^N|\mu,I_t^N,Q_t,\eta_t) \\
    &= P(\mathcal{P}(\X),I^N_{t+1}|\mu,I^N_t,Q_t,\eta_t),
\end{align}
where in the last line we are taking the marginal of the true transition kernel $P(d\kappa' | \kappa,Q,\eta)$ over its first coordinate (the $\pi_{t-N}$ coordinate). We define the cost similarly:
\begin{equation}
     c_N(\hat{\kappa},Q,\eta) = \sum_{(x,u,q)\in\X\times\mathbb{U}\times\M} \eta(q,Q)Q(q|x)\psi(\hat{\kappa})(x)c(x,u).
\end{equation}

Then $\text{MDP}_N \coloneqq (\hat{\W},\mathcal{Q} \times \mathcal{H},P_N(\hat{\kappa}' | \kappa, Q, \eta), c_N(\hat{\kappa},Q,\eta))$. The discounted cost for $\text{MDP}_N$ will be $\hat{J}_\beta(\hat{\kappa}_0,\hat\gamma)$ and the optimal discounted cost  is $\hat{J}_\beta(\hat{\kappa}_0)$. By making $P_N$ and  $c_N$ constant over all of $\cal{P}(\X)$ we will extend the optimal cost to $\tilde{J}_\beta(\kappa_0)$. Let us say that $\hat{\gamma}^*_N$ achieves this optimal cost for the $\text{MDP}_N$ (which certainly exists here the space is finite \cite{yuksel2020control}).

To study the performance of the approximation, we must consider the conditions that $\hat\pi_t$ can be used as an appropriate replacement for $\pi_t$. For a large $N$, this relies on filter stability which measures how quickly a process can recover from starting from the incorrect prior \cite{cregg2024reinforcementlearningoptimaltransmission}. Let us define:
\begin{equation}\label{loss term}
    L_t^N = \sup_{\gamma\in\Gamma_{C-FW}} E^\gamma_\mu [\norm{\pi_t-\hat\pi_t}_{TV}]
\end{equation}
which measures the maximum total variation distance between the predictor with the correct and incorrect prior, at time $t$. From \cite{dobrushin1956central}, we define the Dobrushin coefficient for $P \coloneqq P(x' | x,u)$:
\begin{equation}\label{dobrushing coeff}
    \delta(P,u) = \min_{i,k \in \X}\sum_{j\in\X}\min(P(j|i,u),P(j|k,u)).
\end{equation} 
We can now state the following, extending \cite[Theorem 3.6]{MYDobrushin2020}, and recalling that $\overline{\pi}_t$ is the \emph{filter} obtained by further conditioning $\pi_t$ on $q_t$. Also, to make the dependence on the initialization of $\pi_0$ clear, we write $\pi^\mu_t$ to denote the predictor process when $\pi_0 = \mu$ (and similarly for the filter process).
\begin{lemma}
    \begin{align*}
        & \int_\X \int_\M \norm{\overline\pi_t^\mu-\overline\pi_t^\nu}_{TV}P_\mu^\gamma (dq_t|x_t,q_{[0,t-1]})\pi^\mu_t(dx_t) \\
        & \leq 2\norm{\pi_t^\mu-\pi_t^\nu}_{TV}.
    \end{align*}
\end{lemma}
The proof of this result follows closely that of \cite[Lemma 3]{cregg2024reinforcementlearningoptimaltransmission}. This proof involves taking the Dobrushin coefficient of the conditional probability of the channel output given the state $x$ and the quantizer $Q$. We note that, since the channel is noiseless, the channel kernel $O(q'|q)$, would be an identity map and $\delta(I) = 0$. We note that if the channel were noisy, the right side of this inequality would include a multiplier of $(2-\delta(O))$.
\begin{theorem}\label{loss bound}
For any $\mu\ll\nu$ and $\gamma\in\Gamma_{C-P}$,
    \begin{equation}
        E^\gamma_{\mu} [\norm{\pi_{t+1}-\hat\pi_{t+1}}_{TV}] \leq 2(1-\min_{u\in\mathbb{U}}\delta(P,u))E_\mu^\gamma[\norm{\pi_t-\hat\pi_t}_{TV}]
    \end{equation}
\end{theorem}
\textbf{Proof.}
     To prove this, we will show the statement for two instances of $\pi$ starting from distict priors $\mu$ and $\nu$, where $\mu\ll\nu$. Note that $\overline\pi_t^\mu$ can always be recursively computed when given $\mu$, $q_{[0,t]}$, under some policy $\gamma \in \Gamma_{C-P}$, so it is enough to only take the expectation with respect to $q_{[0,t]}$.
    \begin{align*}
        & E_\mu^\gamma[\norm{\overline\pi_t^\mu-\overline\pi_t^\nu}_{TV}] \\
        &= \int_{\M^t}\int_\M \norm{\bar\pi_t^\mu - \bar\pi_t^\mu}_{TV} P_\mu^\gamma(dq_{[0,t]})\\
        &= \int_{\M^t}\int_\X \int_\M \norm{\bar\pi_t^\mu - \bar\pi_t^\mu}_{TV} P_\mu^\gamma(dq_t|x_t,q_{[0,t-1]}) \\ & \hspace{4cm} P_\mu^\gamma(dx_t|q_{[0,t-1]})P_\mu^\gamma(dq_{[0,t-1]})\\
        &\leq 2\int_{\M^t} \norm{\pi_t^\mu-\pi_t^\nu}P_\mu^\gamma(dq_{[0,t-1]})\\
        &= 2E^\gamma_\mu[\norm{\pi_t^\mu-\pi_t^\nu}]
    \end{align*}
    The second last line is from Lemma \ref{loss bound}. 
    Since the Dobrushin coefficient is a contraction, we can state the following \cite{dobrushin1956central}:
    \begin{align*}
        & E_\mu^\gamma [\norm{\pi_{t+1}^\mu-\pi_{t+1}^\nu}_{TV}]\\
        \leq &(1-\min_{u\in\mathbb{U}}\delta(P,u)E_\mu^\gamma[\norm{\overline{\pi}_t^\mu-\overline{\pi}_t^\mu}_{TV}]\\
        \leq &2(1-\min_{u\in\mathbb{U}}\delta(P,u)E^\gamma_\mu[\norm{\pi^\mu_t-\pi_t^\nu}_{TV}]
    \end{align*} \qed

A bound for $L_t^N$ is now \cite{cregg2024reinforcementlearningoptimaltransmission}:
\begin{align*} 
    L_t^N &\leq \sup_{\gamma\in\Gamma_{C-P}} (2(1-\min_{u\in\mathbb{U}}\delta(P,u)))^N E_\mu^\gamma[\norm{\pi_t-\hat\pi_t}_{TV}]\\
    &\leq (2(1-\min_{u\in\mathbb{U}}\delta(P,u)))^N \norm{\pi_t-\hat\pi_t}_{TV} \\
    &\leq 4(1-\min_{u\in\mathbb{U}}\delta(P,u)))^N
\end{align*}
We will use this loss bound to give a performance bound that shows the impact of using the finite window approximation. The following are extensions of \cite[Theorems 1 and 2]{cregg2024reinforcementlearningoptimaltransmission}. We will define $\norm{c}_\infty \coloneq \max_{(x,u)\in\X\times\mathbb{U}} c(x,u)$.
\begin{theorem} \label{T5.3}
    For any $\gamma\in\Gamma_{C-FW}$ which acts on the first $N$ time steps to generate $\kappa_0$ and any $\pi_{-N}\in \cal{P}(\X)$,
    \begin{equation*}
        E_{\pi_{-N}}^\gamma \big[|\tilde{J}_\beta(\kappa_0)-J_\beta(\kappa_0)]\leq \frac{2\norm{c}_\infty}{(1-\beta)^2}(2(1-\min_{u\in\mathbb{U}}\delta(P,u)))^N
    \end{equation*}
\end{theorem}
\textbf{Proof.}\\
We will show this for $N=1$. Note that $I_t^1 = (q_t,Q_t,\eta_t)$. By definition of $\tilde{J}_\beta$, we have $\tilde{J}_\beta(\kappa_0) = \hat{J}_\beta(\hat{\kappa}_0)$, and thus by the fixed-point equation for $\hat{J}_\beta$ (see e.g., \cite[Chapter 4.2]{HernandezLermaMCP}), we have
\begin{align*}
\sum_{q_1 \in \mathcal{M}}& \tilde{J}_{\beta}(\pi_0, (q_0, Q_0,\eta_0)) P(q_1 | \kappa_0, Q_0,\eta_0)\\
&= \sum_{q_1 \in \mathcal{M}} \hat{J}_{\beta}(\mu, (q_0, Q_0,\eta_0)) P(q_1 | \kappa_0, Q_0,\eta_0).
\end{align*}
We add and subtract the above term, and use the fixed-point equations, to arrive at:
\begin{align*}
&E_{\pi_{-N}}^{\gamma} \left[ \left| J_{\beta}(\kappa_0) - \tilde{J}_{\beta}(\kappa_0) \right| \right] \\
&\hspace{0.2em}\leq ( \norm{c}_{\infty} + \beta \norm{\hat{J}_{\beta}}_{\infty} ) L_0^1
+ \sup_{\gamma'\in\Gamma_{C-FW}} \beta E_{\pi_{-N}}^{\gamma'} \left[ \left| J_{\beta}(\kappa_1) - \tilde{J}_{\beta}(\kappa_1) \right| \right]
\end{align*}
where $L_0^1$ is defined as in \eqref{loss term}. We apply the same process on the supremum term and recursively arrive at:
\begin{align*}
E_{\pi_{-N}}^{\gamma} \left[ \left| J_{\beta}(\kappa_0) - \tilde{J}_{\beta}(\kappa_0) \right| \right] \leq \frac{\norm{c}_\infty}{1-\beta} \sum^\infty_{k=0} \beta^k L_k^1
\end{align*}
where we used the fact that $\norm{J_{\beta}}_{\infty} \leq \frac{\norm{c}_{\infty}}{1-\beta}$. Using the results of Theorem \ref{loss bound} and the geometric series, we have the desired result:
\begin{align*}
E_{\pi_0}^{\gamma} \left[ \left| J_{\beta}(\kappa_0) - \tilde{J}_{\beta}(\kappa_0) \right| \right] &\leq \frac{\norm{c}_\infty}{1-\beta} \sum^\infty_{k=0} \beta^k L^1_k \\
& \le \frac{\norm{c}_\infty}{(1-\beta)^2} 4(1-\min_{u\in\mathbb{U}}\delta(P,u)))^N
\end{align*}
\qed\\
The following theorem gives a performance bound and an exponential convergence rate (in window length) for the approximate model.
\begin{theorem} \label{T14}
    Let $\tilde\gamma_N^*$ be extended from $\hat\gamma^*_N$ over $\cal{P}(\X)$, where $\hat\gamma^*_N$ is optimal for $MDP_N$. Then for any $\gamma \in \Gamma_{C-P}$ which is applied $N$ times starting from $\pi_{-N}$ to generate $\kappa_0$, we have
    \begin{equation*}
        E_{\pi_{t-N}}^\gamma \big[J_\beta(\kappa_0,\tilde\gamma_N^*)-J_\beta(\kappa_0)]\leq \frac{4\norm{c}_\infty}{(1-\beta)^2}(2(1-\min_{u\in\mathbb{U}}\delta(P,u)))^N
    \end{equation*}
\end{theorem}
\textbf{Proof.}\\
First we will state (again for $N=1$):
\begin{align*}
    &E_{\pi_{t-N}}^\gamma \big[J_\beta(\kappa_0,\tilde\gamma_N^*)-J_\beta(\kappa_0)]\\
    &\leq E_{\pi_{t-N}}^\gamma \big[J_\beta(\kappa_0,\tilde\gamma_N^*)-\tilde{J}_\beta(\kappa_0)] + E_{\pi_{t-N}}^\gamma \big[\tilde{J}_\beta(\kappa_0)-J_\beta(\kappa_0)]
\end{align*}
Using a process similar to the one used to prove Theorem \ref{T5.3}, we can get
\begin{align*}
     E_{\pi_{-N}}^\gamma \big[J_\beta(\kappa_0, \tilde{\gamma}^*_N) - \tilde{J}_\beta(\kappa_0)]\leq \frac{2\norm{c}_\infty}{(1-\beta)^2}(2(1 - \min_{u\in\mathbb{U}}\delta(P,u)))^N
\end{align*}
and thus, 
\begin{align*}
    E_{\pi_{t-N}}^\gamma \big[J_\beta(\kappa_0,\tilde\gamma_N^*)-J_\beta(\kappa_0)] \le \frac{4\norm{c}_\infty}{(1-\beta)^2}(2(1-\min_{u\in\mathbb{U}}\delta(P,u)))
\end{align*}
\qed
\subsection{Sliding Finite Window Q-Learning}
Let us now consider the convergence of Q-learning to a solution for the finite sliding window $\text{MDP}_N$. Consider the following sequences:
\begin{enumerate}
    \item $\{\hat{\kappa}_t\}_{t\geq0}$ as defined above.
    \item $\{Q_t\}_{t\geq0}$ is chosen uniformly from $\cal{Q}$ at each time t. 
    \item $\{\eta_t\}_{t\geq0}$ is chosen uniformly from $\cal{H}$ at each time t. 
    \item $\{C_t\}_{t\geq0}$ where $C_t = c_N(\hat{\kappa}_t,Q_t,\eta_t)$
    \item $\{\alpha_t\}_{t\geq0}$, where $\alpha_t: \hat{\W}\times\cal{Q}\times\cal{H}\mapsto \R_+$ is the learning rate at time $t$. Specifically, we define $$\alpha_t(\hat{\kappa},Q,\eta) = \frac{1}{1+\sum^t_{k=0} \mathbf{1}_{\{(\hat{\kappa}_k,Q_k,\eta_k) = (\hat{\kappa},Q,\eta)\}}}$$
    if $(\hat{\kappa}_t, Q_t, \eta_t) = (\hat{\kappa}, Q, \eta)$, and $0$ otherwise.
    \item $\{\mathbf{Q}_t\}_{t\geq0}$, where $\mathbf{Q}_t: \hat\W\times\cal{Q}\times\cal{H}\mapsto \R_+$, $\mathbf{Q}_0 \equiv 0$, and $\mathbf{Q}_t$ is updated as:
    \begin{align*}
    & \mathbf{Q}_{t+1}(\hat{\kappa},Q,\eta)
    = (1-\alpha_t(\hat{\kappa},Q,\eta))\mathbf{Q}_t(\hat{\kappa},Q,\eta)\\
    & \hspace{1cm} +\alpha_t(\hat{\kappa},Q,\eta)[C_t +\beta \min_{(Q,\eta)\in\cal{Q}\times\cal{H}} \mathbf{Q_t}(\hat{\kappa}_{t+1},Q,\eta)] 
\end{align*} 
\end{enumerate}
We will argue that the sequences $(\hat{\kappa}_t,Q_t,\eta_t,C_t,\alpha_t)$ satisfy Assumption \ref{Qconvergence assumption}, and thus the Q-Learning iterations will converge to a meaningful limit. As in the previous section, we require Assumption \ref{assumption:PHR}, which we recall is that $\{x_t\}_{t \ge 0}$ is positive Harris recurrent (and thus has unique invariant measure $\zeta$) under the uniform choice of $(Q_t, \eta_t)$.

In the following result, when we say \emph{almost any} $(\hat{\kappa},Q,\eta)\in \hat{\W}\times\mathcal{Q}\times{\mathcal{H}}$, we mean any $(\hat{\kappa},Q,\eta)$ which has positive probability of occurring under any prior; under Assumption \ref{assumption:PHR}, this is equivalent to it having positive probability under the unique invariant measure $\zeta$.

\begin{lemma} For any $\pi_{-N}$ and almost any $(\hat{\kappa},Q,\eta)\in \hat{\W}\times\cal{Q}\times{\mathcal{H}}$, under Assumption \ref{assumption:PHR} we have that almost surely:
\begin{enumerate}[(i)]
    \item $(\hat{\kappa}_t,Q_t,\eta_t) = (\hat{\kappa},Q,\eta)$ infinitely often and $\sum_{t\geq0} \alpha_t(\hat{\kappa},Q,\eta)=\infty$.
    \item 
        $$\frac{\sum^t_{k=0} C_k \mathds{1}_{\{\hat{\kappa}_k=\hat{\kappa},Q_k=Q,\eta_k=\eta\}}}{\sum^t_{k=0}  \mathds{1}_{\{\hat{\kappa}_k=\hat{\kappa},Q_k=Q,\eta_k=\eta\}}} \rightarrow c_N(\hat{\kappa},Q,\eta)$$ 
    \item For any $g$,
        \begin{multline*}
            \frac{\sum^t_{k=0} g(\hat{\kappa}_{k+1}) \mathds{1}_{\{\hat{\kappa}_k=\hat{\kappa},Q_k=Q,\eta_k=\eta\}}}{\sum^t_{k=0} \mathds{1}_{\{\hat{\kappa}_k=\hat{\kappa},Q_k=Q,\eta_k=\eta\}}} \\ \rightarrow \int_{\hat{\W}} g(\hat{\kappa}') P_N(d\hat{\kappa}'|\hat{\kappa},Q,\eta)
        \end{multline*}
\end{enumerate}
\end{lemma}
\textbf{Proof.}
    Under Assumption \ref{assumption:PHR}, we have the following as $t\rightarrow \infty$:
    $\norm{P(x_t\in\cdot) - \zeta}_{TV}\rightarrow 0$
    as $t\rightarrow \infty$. Since for almost every $\hat{\kappa}$ we have that $P_\zeta(\hat{\kappa}) > 0$, we will eventually have $P_{\pi_{t-N}}(\hat{\kappa}) > 0$ and thus almost every $(\hat{\kappa},Q,\eta)$ is hit infinitely often, satisfying (i). Part (ii) is immediate because $C_t = c_N(\hat{\kappa},Q,\eta)$. 
    
    On part (iii); since the marginals of $x_t$ converge to $\zeta$, under the exploration policy, the marginals on $I_t^N$ also converge and we have:
    \begin{align*}
    &\frac{\sum^t_{k=0} g(\hat{\kappa}_{k+1}) \mathds{1}_{\{\hat{\kappa}_k=\hat{\kappa},Q_k=Q,\eta_k=\eta\}}}{\sum^t_{k=1} \mathds{1}_{\{\hat{\kappa}_k=\hat{\kappa},Q_k=Q,\eta_k=\eta\}}} \\
    &= \frac{\sum^t_{k=0} g(I_{k+1}) \mathds{1}_{\{I_k=i,Q_k=Q,\eta_k=\eta\}}}{\sum^t_{k=1} \mathds{1}_{\{I_k=i,Q_k=Q,\eta_k=\eta\}}} \\
    &\rightarrow \int g(i')P(\pi,i'|\zeta,i,Q,\eta) \\
    & = \int_{\hat{\W}} g(\hat{\kappa}') P_N(d\hat{\kappa}'|\hat{\kappa},Q,\eta)
\end{align*}\qed

We can now state the following as a result of Theorem \ref{Q-learning convergence theorem}. 
\begin{theorem}\label{theorem:finite_memory_Qlearning}
    If $\delta(P,u)>1/2$ for all $u \in \mathbb{U}$, and under Assumption \ref{assumption:PHR}, we have the following: for almost every $(\hat{\kappa},Q,\eta) \in \hat{\W}\times \cal{Q} \times\cal{H}$, $\mathbf{Q}_t(\hat{\kappa},Q,\eta)$ will converge almost surely to a limit satisfying: 
    \begin{align*}
        \mathbf{Q}^*(&\hat{\kappa},Q,\eta) = c_N(\hat{\kappa},Q,\eta)\\&+\beta\sum_{\hat{\kappa}'\in \hat{\W}} P_N(\hat{\kappa}'\mid \hat{\kappa},Q,\eta) \min_{(Q',\eta')\in(\cal{Q}\times\cal{H})} \mathbf{Q}^*(\hat{\kappa'},Q',\eta')
    \end{align*}
    This coincides with the DCOE for $\text{MDP}_N$, and therefore the policy $\hat{\gamma}^*_N(\hat{\kappa})= \argmin_{Q,\eta\in{\cal{Q}\times\cal{H}}}\mathbf{Q}^*(\hat{\kappa},Q,\eta)$ is optimal for $\text{MDP}_N$. Thus if we define a policy $\tilde{\gamma}_N^*(\kappa) \coloneq \hat{\gamma}^*_N(\hat{\kappa})$, we have by Theorem \ref{T14} that for almost every $\kappa_0$,
$$ \lim_{N\to\infty}|J_\beta(\kappa_0,\tilde{\gamma}_N^*) - J_\beta(\kappa_0)|=0.$$
\end{theorem}

\section{Comparison of Finite Predictor State and Finite Window Approximations}\label{comparisonModelsSec}
We will now compare the two rigorously justified approximation methods presented: The finite sliding window method and quantized predictor approach.


\begin{itemize}
\item[(i)] [Computational efficiency] The finite sliding window method offers the benefit that all possible values of $\hat{\kappa}_t$ and $\hat{\pi}_t=\psi(\hat{\kappa}_t)$ can be computed offline before running the learning algorithm. 
However, the space of possible windows grows exponentially in the underlying state or action space dimension and the window length parameter. Because of this, each iteration requires less computation but the number of Q-iterations before the algorithm converges to an optimal policy can be very large for window sizes larger than 2. Conversely, the computational complexity at each Q-iteration in the implementation of the quantized predictor approach is significantly higher. Learning requires frequent updates of the predictor state (and requires a Bayesian update which necessitates access to the system model), and as the quantization becomes finer, the number of possible states grows exponentially (although not every state may be relevant, as we only need consider those with positive measure under the unique invariant measure, see Theorem \ref{theorem:quantizedQlearning_nearoptimal}). Thus, even for resolutions of quantization up to 15, the number of visited states is small enough that the algorithm can converge faster than the sliding window method (for window sizes greater than 2). This leads to a trade-off between approximation accuracy and computational costs for both these methods.
\item[(ii)] [Universality, and availability of model information] The predictor quantization method requires Bayesian updates of the predictor state (which necessitates access to the system model). The sliding window method does not need access to the system model and is thus universal as long as the required mixing conditions are known to hold.

\item[(iii)] [Initialization] The sliding window method is insensitive to initialization. Provided the window length is sufficiently large, it converges to near-optimal policies regardless of the initial distribution, making it advantageous in real-world applications where the initial state might not be precisely known or controllable. On the other hand, the quantized predictor state space method is sensitive to initialization. The learned policies from this approach are near-optimal only when the initial state has positive measure under the unique invariant measure for $\pi_t$. This may be problematic in practical scenarios where the system may start from an arbitrary initial condition.

\item[(iv)] [On required predictor/filter stability properties] The finite sliding window method assumes a strong form of controlled predictor/filter stability. Only under appropriate conditions, such as those characterized by the Dobrushin coefficient, can this approach have near-optimal performance. The quantized state space method works under weaker conditions of stability, and does not require an assumption on the Dobrushin coefficient.
\end{itemize}


Thus, both methods achieve near-optimality under the right conditions, but they are suitable in complementary settings: The sliding-window approach can be more computationally efficient and more robust to changes in initial conditions, making it preferable for real-time or data-driven applications—such as autonomous vehicle or drone coordination—where model information is limited and faster operation is important. The predictor-quantization method, while computationally less efficient, is flexible with regard to filter stability conditions and useful for model-based, low-dimensional systems such as industrial process control or robotics, where system dynamics are relatively better known.

\section{Simulations}

In this section, we provide simulation studies where both of the approximation methods introduced in the paper are implemented in the following.

\subsection{Simulation: Finite Model Approximation via Predictor Quantization}\label{quantizedbeliefsimulations}
We will now give an example of the control problem and simulate the performance of the algorithm using the finite state approximation. We will use a discount factor of $\beta = 0.8$.

Let $\X = \{1,2,3\}$, $\mathbb{U} = \{1,2\}$, and $\M = \{1,2\}$. The transition kernel $P(\cdot|x,u)$ is as follows:
\begin{align*}
    P(\cdot|x,1) &=\frac{1}{10}
    \begin{pmatrix}
        4&0&6\\
        3&7&0\\
        2.5&2.5&5\\
    \end{pmatrix}\\
    P(\cdot|x,2) &= \frac{1}{10}
    \begin{pmatrix}
        3.5&5&1.5\\
        5&5&0\\
        0&7.5&2.5\\
    \end{pmatrix}
\end{align*}
The cost function is defined as $c(x,u) = C_{(xu)}$, that is the element at $[x,u]$ of the cost matrix:
\begin{equation}
    C = 
    \begin{bmatrix}
        0&0\\
        0&1\\
        1&1\\
    \end{bmatrix}
\end{equation}

\begin{figure}
    \centering
    \includegraphics[width=0.9\linewidth]{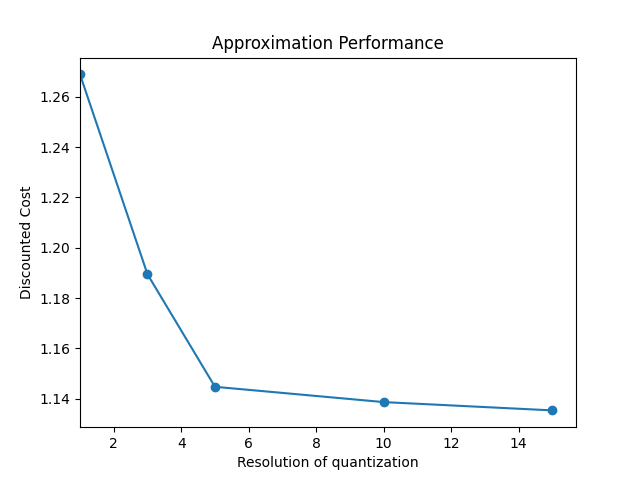}
    \caption{Discounted cost using a learned policy for $\text{MDP}_n$ for quantization resolutions $n=\{1,3,5,10,15\}$}
    \label{fig:beliefquantizationsim}
\end{figure}

By design, the kernel is aperiodic and irreducible (under any stationary control policy for a corresonding fully observable MDP). For a randomly and independently generated first action, $u_0$, the unique invariant distribution of $P(\cdot|x,u_0)$ is computed and used as the initial distribution for learning the best policy. For estimating the cost under this policy, $\pi_0$ is a predictor of the form in \eqref{recurrent_pi}. For quantization resolutions $n=\{1,...5\}$, the aforementioned Q-learning algorithm is applied and policies (after $\sim 10^6$ iterations or earlier convergence) are applied. We calculate the empirical expected discounted cost for the learned policy by running the simulation with a horizon $N =1000$ and averaging this cost over 1000 Monte-Carlo iterations. As the resolution of quantization increases, the cost decreases, as expected. The results are shown in Figure \ref{fig:beliefquantizationsim}. 

\subsection{Simulation: Finite Model Approximation via the Sliding Window Method}
\begin{figure}
    \centering
    \includegraphics[width=0.9\linewidth]{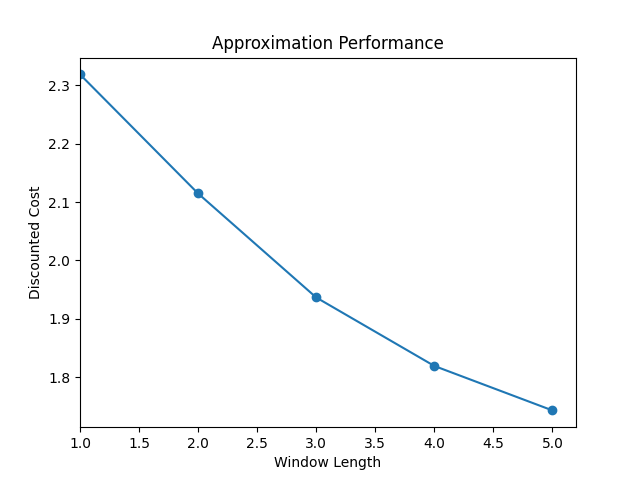}
    \caption{Discounted cost using a learned policy for $MDP_N$ for window lengths $N=\{1,2,3,4,5\}$.}
    \label{fig:SFWResults}
\end{figure}

To simulate the performance of the algorithm using the finite sliding window approximation, we will use the same $\beta$, $\X$, $\mathbb{U}$, $\M$, and $c(x,u)$, as in Section \ref{quantizedbeliefsimulations} above. This example will use the following kernel.
\begin{align}
  P(\cdot|x,1) &=\frac{1}{10}
    \begin{pmatrix}
        4&0&6\\
        4.5&3&2.5\\
        2.5&2.5&5\\
    \end{pmatrix} \nonumber \\
    P(\cdot|x,2) &=\frac{1}{10} 
    \begin{pmatrix}
        3.5&5&1.5\\
        5&5&0\\
        0.5&7.5&2\\
    \end{pmatrix} \nonumber
\end{align}
This kernel satisfies the irreducibility and aperiodicity condition given in Assumption \ref{assumption:PHR} but also satisfies the Dobrushin coefficient condition in Theorem \ref{T5.3}: \[\min_{u\in\mathbb{U}} \delta(P,u) = 0.55\ge 1/2\] which means the results of Theorems \ref{loss bound} and \ref{T5.3} are applicable.\\
The source begins from $x_{-N} = 1$, and the first $N$ steps are obtained by using a uniform choice of $(Q_t, \eta_t)$. A learned policy for the $\text{MDP}_N$ for window lengths $N=\{1,2,3,4,5\}$ was found by training the Q-learning algorithm for $10^5$ iterations. Using this policy, the empirical expected cost was computed using $N = 1000$ and 1000 Monte-Carlo iterations. As the window length increases, the discounted cost decreases consistently, as expected. The results are shown in Figure \ref{fig:SFWResults}.

\section{Extensions}
In this section, we comment on two generalizations of our work presented in the paper. While we do not study them in detail due to space constraints, we provide a detailed discussion on how the extensions can be made with little conceptual effort, though with some additional technical analysis as explicitly detailed in the following. 

\subsection{Case with Noisy Channels and Feedback}\label{noisyChannelExt}

\begin{figure}[h]
        \begin{center}
\epsfig{figure=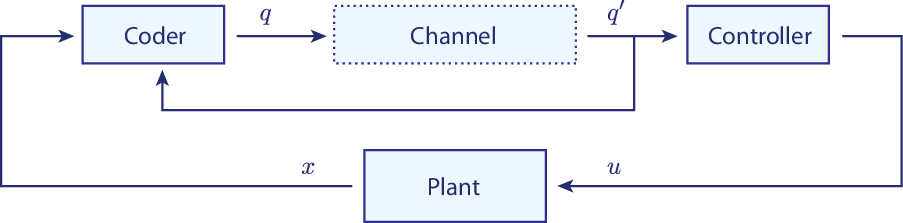,width=9cm}
\caption{Optimal coding and control over a noisy channel. \label{LLL0001}}
        \end{center}
\end{figure}

A more common occurrence in networked control applications involves a noisy channel: Consider the case in Figure \ref{LLL0001} in which there is a noisy channel between the quantization output $q_t$ and the channel output $q'_t$ where the channel is a Discrete Memoryless Channel (DMC), which satisfies the property that $P(q'_t | q_{[0,t]}) = P(q'_t | q_t)$ for all realizations and history. 

In this case, a coding policy $\gamma^{e} = \{\gamma^e_t, t \ge 0\}$ is a sequence of functions which generate quantization outputs, $q_t$, as a measurable function of the encoder's information at time $t$:
\begin{equation}
    I_t^e = \{{x_{[0,t]},q_{[0,t-1]},q'_{[0,t-1]}}\}
\end{equation}
\begin{equation}
    \gamma^e_t:I_t^e\mapsto q_t \in \M
\end{equation}
where $\M$ is a finite quantization output alphabet, so that $q_t = \gamma^e_t(I_t^e)$ for $t \in \mathbb{Z}_+$.


A control policy $\gamma^c = \{\gamma^c_t,t \ge 0\}$ is a sequence of functions which generate control actions, $u_t$, as a measurable function of the information available at the controller at time $t$:
\begin{equation}
    I_t^c = \{{q'_{[0,t]}, u_{[0,t-1]}}\}
\end{equation}
\begin{equation}
    \gamma^c_t:I_t^c \mapsto u_t \in \mathbb{U} 
\end{equation}
where $\mathbb{U}$ is the finite action space, so that $u_t = \gamma^c_t(I_t^c)$ for all $t \in \mathbb{Z}_+$.

For this setup, under the criterion (\ref{objective})-(\ref{infinitediscountcost}), the controlled separation results of the type given in Theorem \ref{Walrand and Varaiya} hold in the noisy channel case as well by \cite{WalrandVaraiya} and \cite[Theorem 15.3.8]{YukselBasarBook24}, provided that the encoder has access to the realizations $q'_t$ in a causal fashion. 

Accordingly, for such a setup, the structural, existence, and approximation results presented in this paper generalize with essentially identical arguments presented in the paper provided that there is feedback from the channel output to the encoder input, see the analysis \cite[Section IV]{WoodLinderYuksel} and \cite{cregg2024reinforcementlearningoptimaltransmission} for the control-free setup. Therefore, the counterparts of the results on existence, structure and near optimality in Sections \ref{section:MDP_Formulation}-\ref{section:Finite_Window} apply nearly identically. 

One technical detail (for the learning results in Section \ref{section:Predictor_Quantization}) is that in the noisy channel case we do not necessarily have recurrent values of $\pi_t$ (of the form in \eqref{recurrent_pi}). Thus, additional filter stability conditions are needed to ensure the uniqueness of an invariant measure of $\pi_t$ for the exploration process. Another  difference is that the Dobrushin coefficient term in Theorem \ref{T5.3} (and the related supporting results) have to be modified to include the effect of the channel.

To summarize, we have the following remark.
\begin{remark} The following counterparts of the results hold for the noisy channel setup:
    \begin{itemize}
        \item [(i)] For the minimization of the infinite horizion discounted cost \eqref{infinitediscountcost} in the noisy channel case, quantizing the resulting MDP leads to a near-optimal policy, i.e., a noisy channel analog of Theorem \ref{special case quantization performance} holds.
        \item [(ii)] One can obtain an optimal policy for the quantized MDP in (i) through Q-learning, i.e., a noisy channel analog of Theorem \ref{theorem:quantizedQlearning_nearoptimal} holds.
        \item [(iii)] For the minimization of the infinite horizion discounted cost \eqref{infinitediscountcost} in the noisy channel case, using a finite memory of past channel outputs and controls leads to a near-optimal policy, i.e., a noisy channel analog of Theorem \ref{T5.3} holds.
        \item [(iv)] One can obtain an optimal policy for the finite memory MDP in (iii) through Q-learning, i.e., a noisy channel analog of Theorem \ref{theorem:finite_memory_Qlearning} holds.
    \end{itemize}
\end{remark}

\subsection{Case with Real or Polish State $\mathbb{X}$ and Action $\mathbb{U}$ Spaces}\label{contSpaceSec}

Many engineering systems involve continuous vector spaces. While we leave a detailed study of this direction for future work, also in view of space constraints, in the following we highlight the program that the extension would entail.

The structural result we have developed applies also to the case with $\mathbb{X}=\mathbb{R}^d$ and $\mathbb{U}$ a standard Borel control action space; as studied in \cite[Chapter 15]{YukselBasarBook24} building on \cite{YukLinZeroDelay}, under the following assumption.

	\begin{assumption}\label{assumption1} Let a controlled system be expressed in the stochastic realization form
		\begin{equation}\label{dynmc}
		x_{t+1}=f(x_t,u_t,w_t), \quad t=0,1,2,. . .,
		\end{equation}
		where $f : \mathbb{R}^d \times \mathbb{U} \times \mathbb{R}^d \rightarrow \mathbb{R}^d$ is a Borel measurable function and ${w_t}$
		is an independent and identically distributed (i.i.d.) vector
		noise sequence which is independent of $x_0$. It is assumed that
		for each fixed $x, u \in \mathbb{R}^d \times \mathbb{U}$, the distribution of $f(x, u, w_t)$ admits
		the (conditional) density function $\phi(\cdot|x,u)$ (with respect to the
		$d$-dimensional Lebesgue measure) which is positive everywhere.
		Furthermore, $\phi(\cdot|x,u)$ is bounded and Lipschitz uniformly in $x,u$.
	\end{assumption}

		Let $\mathcal{G}$ denote the set of all probability measures on $\mathbb{R}^d$ admitting densities that are bounded by a constant $C$ and Lipschitz with constant $C_1$. Then, viewed as a class of probability density functions, $\mathcal{G}$ is uniformly bounded. In \cite[Lemma 3]{YukLinZeroDelay} it is shown that $\mathcal{G}$ is closed in $\mathcal{P}(\mathbb{R}^d)$ under weak convergence. 

As discussed in \cite{YukLinZeroDelay}, a quantizer
$Q$ with cells $\{B_1,\ldots,B_M\}$  can be
characterized as a stochastic kernel $Q$ from $\mathbb{R}^n$ to
$\{1,\ldots,M\}$ defined by
\[
Q(i|x)= 1_{\{x \in B_i\}}, \quad  i=1,\ldots,M.
\]
We can endow the quantizers with the topology induced by a
stochastic kernel interpretation under the Young topology \cite{YukLinZeroDelay}. If $P$ is a probability measure on
$\mathbb{R}^n$ and $Q$ is a stochastic kernel from $\mathbb{R}^n$ to
${\cal M}$, then $PQ$ denotes the resulting joint probability measure on
$\mathbb{R}^n \times {\cal M}$. That is, a quantizer sequence $Q_n$ converges to $Q$ weakly at $P$ ($Q_n \to Q$ weakly
 at $P$) if $PQ^n \to PQ$ weakly. Similarly,   $Q_n$ converges to $Q$
in total variation at $P$  ($Q_n \to Q$ at $P$ in total variation)  if $PQ^n \to PQ$ in total variation.

To facilitate compactness properties, \cite{GyorgyLinder2,GyorgyLinder} \cite{YukselOptimizationofChannels} restrict the set of quantizers considered by only allowing quantizers having convex quantization bins (cells) $B_i$, $i=1,\ldots,M$.

\begin{assumption}\label{A113}
The quantizers have convex codecells with at most a given number of cells; that is the quantizers live in $\mathcal{Q}_c(M)$, the collection of $k$-cell quantizers with convex cells where $1 \leq k \leq M$.
\end{assumption}

In this context, let $\Gamma_{C-P}$  denote the set of all predictor structured controlled Markov policies which in addition satisfy the condition that all quantizers $Q_t$, $t\ge 0$ have convex cells (i.e., $Q_t\in \mathcal{Q}_c(M)$ for all $t\ge 0$).

For this problem, the topological construction on the joint encoder and controller spaces is slightly more involved as the space of quantizers is no longer finite. Notably, the formulation in $(\ref{pi_to_Q}-\ref{pi_to_eta})$ involves a state-dependent action space. Accordingly, following a measurable selection theorem, one can obtain structural optimality results as in Theorem \ref{Markov property} and then the necessary weak Feller regularity on the resulting MDP kernel $P(d\pi_{t+1} | \pi_t, Q_t, \eta_t)$, similarly to \cite[Lemma 11]{YukLinZeroDelay}. Therefore, an existence result and a counterpart to Theorem \ref{costeqiv and blackwell}, can be obtained.

To obtain finite approximations, using the weak Feller regularity, one can quantize the source space and obtain an approximate Markov model \cite[Section 3]{kara2024near}, show that the approximate model is close to the original model under weak convergence, and therefore show that the solution of the approximate model is near optimal for the original model \cite{KSYContQLearning}. 

Accordingly, the analysis in the current paper, by applying an optimal zero-delay coding and control for an approximate finite model, will be near-optimal for the original model under mild technical conditions. Furthermore, analogous reinforcement learning results can be obtained. 

\begin{remark} [A Separated Design] An alternative approach is to directly approximate a continuous space MDP with a finite space MDP with uniform guarantees on the approximation error under arbitrary piece-wise constant (across bins) policies (so that the approximation error vanishes as the size of the finite model approximation grows) building on the proof of \cite[Theorem 3]{KSYContQLearning} (and the analysis in \cite{SaYuLi15c}), and then to apply the analysis of our paper as if the approximate finite model is the true model. Such a separated design (of model approximation \cite{KSYContQLearning} and then optimal coding and control for the approximate model as established in our paper) leads to a modular approach with guaranteed performance bounds and with only modest additional mathematical analysis, even though it would be suboptimal compared with the joint optimization outlined above (where the coding and control is jointly optimized without an approximation step). 
\end{remark}

In summary, while the construction and proof program are outlined as above, a detailed technical analysis of this practically significant setup is left for future work, with refinements on the technical conditions to be presented and the noted measurable selection theorem above to be precisely stated and proven.  


\section{Conclusion}
This work has addressed the problem of optimal control over a finite-rate noiseless communication channel, proving structural, existence, and approximation results that characterize optimal and near-optimal coding and control policies for a controlled Markovian system. For the purpose of Q-learning, we introduced two methods with guaranteed performance bounds: the finite sliding window and the quantized predictor measurate methods. We demonstrated the effectiveness of Q-learning for both methods and noted the trade-offs between computational efficiency and initialization sensitivity. While the sliding window approach is more practical for real-time systems due to its lower complexity and universality, the quantized predictor measure method allows for finer approximations but at a higher computational cost. Future research problems involve extending the analysis and the methods to more general scenarios, such as noisy channels with and without feedback, and systems with uncountable state and action spaces, as outlined in the paper. 

\bibliographystyle{IEEEtran}


\begin{IEEEbiography}[{\includegraphics[width=1in,height=1.25in,clip,keepaspectratio]{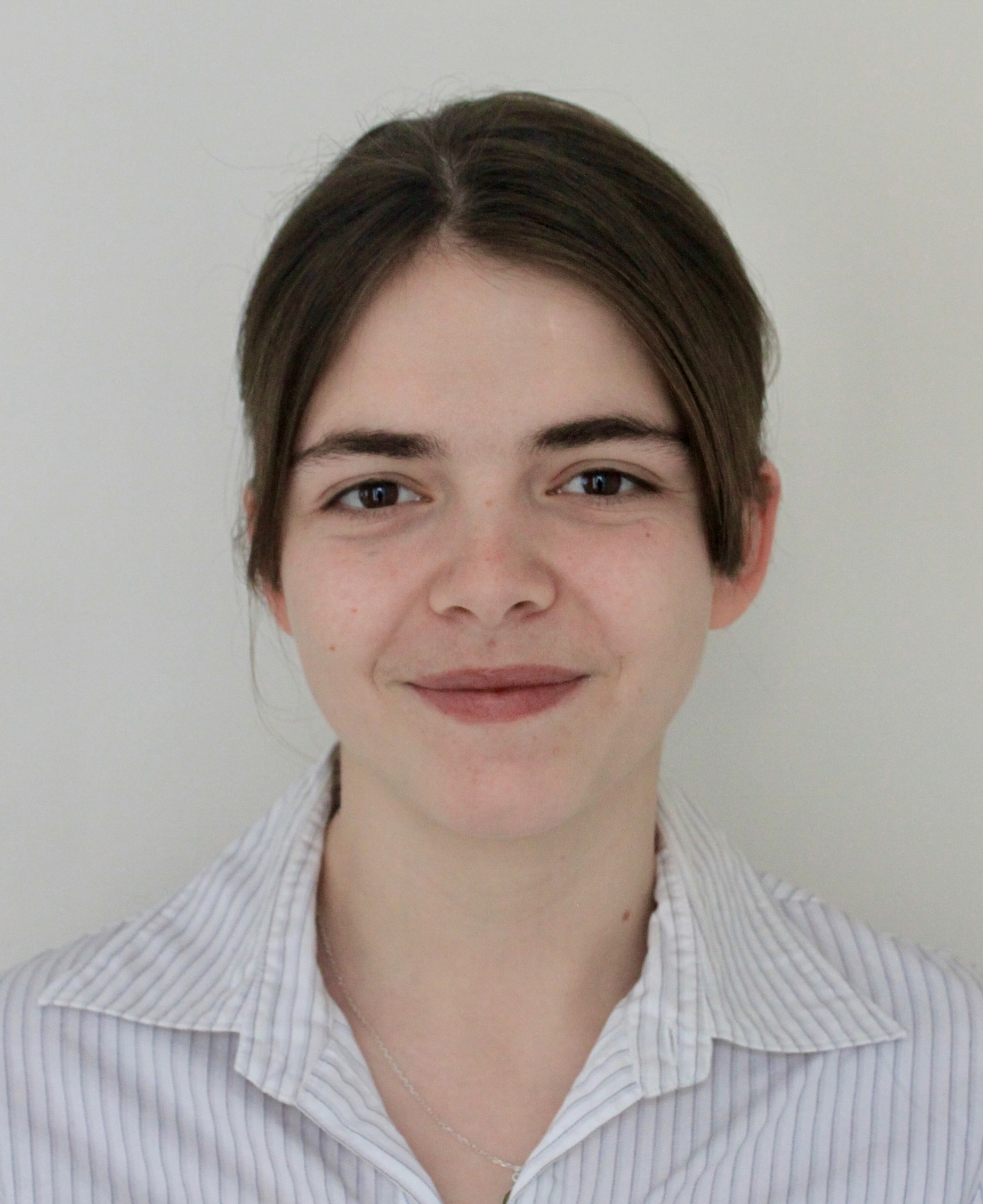}}]
{\bf Evelyn Hubbard} received her B.A.Sc. degree in Mathematics and Engineering. She is currently a Master's student at McGill University in the Department of Mathematics and Statistics. Her research interests include mathematical statistics, stochastic control theory, information theory and probability.
\end{IEEEbiography}

\begin{IEEEbiography}[{\includegraphics[width=1in,height=1.25in,clip,keepaspectratio]{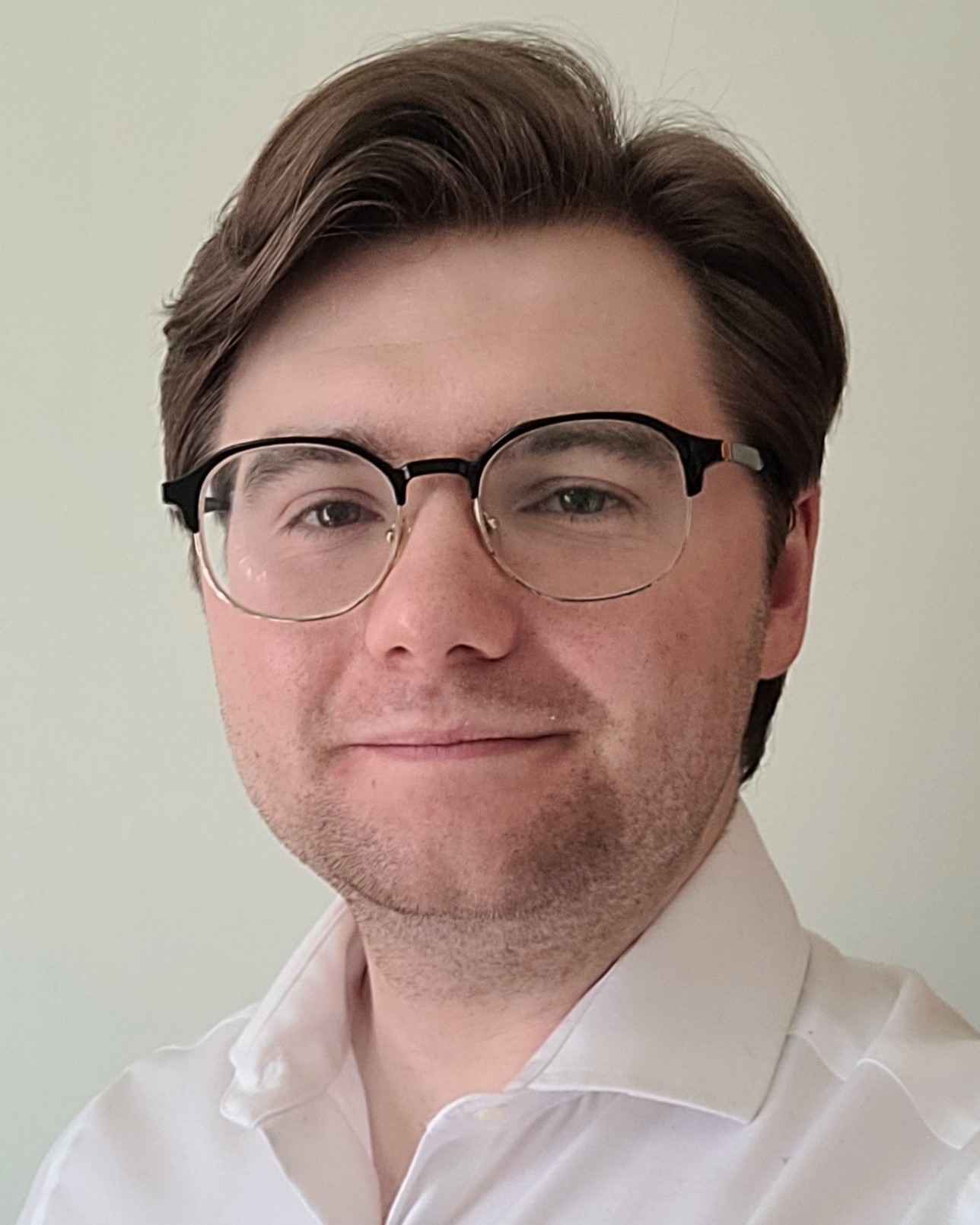}}]
{\bf Liam Cregg} received his B.A.Sc. degree in Mathematics and Engineering and Master's degree in Mathematics and Statistics from Queen's University. He is currently a PhD student with the Institut für Automatik at ETH-Zurich. His research interests include stochastic control theory, information theory and probability.
\end{IEEEbiography}

\begin{IEEEbiography}[{\includegraphics[width=1in,height=1.25in,clip,keepaspectratio]{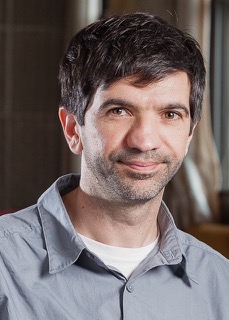}}]
{\bf Serdar Y\"uksel} (S'02, M'11, SM'23) received his B.Sc. degree in Electrical and Electronics Engineering from Bilkent University in 2001; M.S. and Ph.D. degrees in Electrical and Computer Engineering from the University of Illinois at Urbana-Champaign in 2003 and 2006, respectively. He was a post-doctoral researcher at Yale University before joining Queen's University as an Assistant Professor in the Department of Mathematics and Statistics, where he is now a Professor. His research interests are on stochastic control theory, information theory, and probability. 
\end{IEEEbiography}

\end{document}